\documentclass[11pt]{amsart}
\usepackage{amssymb,a4wide}
\usepackage[all]{xy}
\usepackage{graphicx}

\title[Globular realization and cubical underlying homotopy
type]{Globular realization and cubical underlying homotopy type of
  time flow of process algebra}
\author[P. Gaucher]{Philippe Gaucher}
\address{Laboratoire PPS  (CNRS UMR 7126)\\ Universit{\'e} Paris 7--Denis Diderot\\
  Case 7014\\ 75205 PARIS Cedex 13 \\ France}
\email{gaucher@pps.jussieu.fr}
\urladdr{http://www.pps.jussieu.fr/{\~{}}gaucher/}  

\newcommand{\bea}{\begin{eqnarray}}
\newcommand{\eea}{\end{eqnarray}}
\newcommand{\beas}{\begin{eqnarray*}}
\newcommand{\eeas}{\end{eqnarray*}}

\swapnumbers

\newtheorem{thm}{Theorem}[subsection]
\newtheorem*{thmN}{Theorem}
\newtheorem{prop}[thm]{Proposition}
\newtheorem{lem}[thm]{Lemma}
\newtheorem{cor}[thm]{Corollary}

\newtheorem{defn}[thm]{Definition}

\newtheorem{nota}[thm]{Notation}
\newcommand{\bd}{\begin{defn}}
\newcommand{\ed}{\end{defn}}
\newcommand{\bp}{\begin{prop}}
\newcommand{\ep}{\end{prop}}
\newcommand{\bth}{\begin{thm}}
\renewcommand{\eth}{\end{thm}}
\newcommand{\bpf}{\begin{proof}}
\newcommand{\epf}{\end{proof}}

\newcommand{\fl}[1]{\ar@{->}[l]_-{#1}}
\newcommand{\fr}[1]{\ar@{->}[r]^-{#1}}
\newcommand{\fd}[1]{\ar@{->}[d]_-{#1}}
\newcommand{\fu}[1]{\ar@{->}[u]^-{#1}}

\newcommand{\C}{\mathcal{C}}
\newcommand{\D}{\mathcal{D}}
\newcommand{\N}{\mathbb{N}}
\newcommand{\de}{\partial}
\newcommand{\p}\times
\renewcommand{\vec}{\overrightarrow}
\renewcommand{\P}{\mathbb{P}}
\newcommand{\brm}[1]{\rm{\mathbf{#1}}}
\newcommand{\ho}{{\mathbf{Ho}}}

\newcommand{\iso}{\cong}
\newcommand{\lp}{\left(}
\newcommand{\rp}{\right)}
\newcommand{\vI}{\vec{I}}

\renewcommand{\leq}{\leqslant}
\renewcommand{\geq}{\geqslant}

\renewcommand{\top}{{\brm{Top}}}
\newcommand{\gltop}{{\brm{glTop}}}
\newcommand{\dtop}{{\brm{Flow}}}
\newcommand{\set}{{\brm{Set}}}
\newcommand{\tdtop}{{\brm{FLOW}}}
\newcommand{\ttop}{{\brm{TOP}}}
\newcommand{\tgltop}{{\brm{glTOP}}}
\newcommand{\mtop}{{\brm{MTop}}}
\newcommand{\glob}{{\rm{Glob}}}
\newcommand{\ot}{\otimes}
\newcommand{\sis}{\Delta^{op}\set}
\DeclareMathOperator{\id}{Id}
\DeclareMathOperator{\diag}{{\rm{Diag}}}
\DeclareMathOperator{\cell}{{\brm{cell}}}

\newcommand{\hda}{{\cell(I^{gl}_+)}}
\newcommand{\ddownarrow}{{\downarrow}}
\DeclareMathOperator{\sing}{Sing}

\DeclareMathOperator{\gl}{{\rm{gl}}}
\DeclareMathOperator{\hogl}{{\rm{hogl}}}
\DeclareMathOperator{\map}{{\rm{Map}}}

\def\cartesien{%
  \ar@{-}[]+R+<6pt,-2pt>;[]+RD+<6pt,-6pt>%
  \ar@{-}[]+D+<2pt,-6pt>;[]+RD+<6pt,-6pt>%
}
\def\cocartesien{%
  \ar@{-}[]+L+<-6pt,+2pt>;[]+LU+<-6pt,+6pt>%
  \ar@{-}[]+U+<-2pt,+6pt>;[]+LU+<-6pt,+6pt>%
}
\def\hocartesien{%
  \ar@{-}[]+R+<6pt,-2pt>;[]+RD+<6pt,-6pt>_{h}%
  \ar@{-}[]+D+<2pt,-6pt>;[]+RD+<6pt,-6pt>%
}
\def\hococartesien{%
  \ar@{-}[]+L+<-6pt,+2pt>;[]+LU+<-6pt,+6pt>_{h}%
  \ar@{-}[]+U+<-2pt,+6pt>;[]+LU+<-6pt,+6pt>%
}

\newcommand{\liminj}{\varinjlim}
\newcommand{\limproj}{\varprojlim}

\begin{document}

\begin{abstract}
  We construct a small realization as flow of every precubical set
  (modeling for example a process algebra). The realization is small
  in the sense that the construction does not make use of any
  cofibrant replacement functor and of any transfinite construction.
  In particular, if the precubical set is finite, then the
  corresponding flow has a finite globular decomposition. Two
  applications are given.  The first one presents a realization
  functor from precubical sets to globular complexes which is
  characterized up to a natural S-homotopy.  The second one proves
  that, for such flows, the underlying homotopy type is naturally
  isomorphic to the homotopy type of the standard cubical complex
  associated with the precubical set.
\end{abstract}

\maketitle

\tableofcontents

\section{Introduction}

\subsection{Presentation of the results} 
 
Various topological models of concurrency \cite{survol} have been
introduced so far.  Local pospaces \cite{MR1683333} are topological
spaces equipped with a local partial ordering representing a time
flow.  $D$-spaces \cite{mg} are topological spaces equipped with a
family of continuous paths playing the role of execution paths.  A
$d$-space is not necessarily locally partially ordered but their
category is complete and cocomplete. This is an advantage of this
model on the one of local pospaces. A close framework is the full
subcategory of $d$-spaces which are colimit-generated by a small full
subcategory of cubes \cite{FR}. The interest of the latter category is
that it is locally presentable, and that it is therefore possible to
construct directed coverings using strict factorization system
techniques. Streams \cite{SK} are locally preordered topological
spaces.  The category of streams is also complete and cocomplete.
Every $d$-space and every local pospace can be viewed as a stream.
Finally, the globular complexes \cite{diCW} \cite{model2} are
topological spaces equipped with a globular decomposition which is the
directed analogue of the cellular decomposition of a CW-complex. The
globular complexes can be viewed as a subcategory of the categories of
local po-spaces, of $d$-spaces, of $d$-spaces colimit-generated by
cubes, and of streams.  Nevertheless, the category of globular
complexes is big enough to contain all examples coming from
concurrency \cite{diCW}. All these models start from a topological
space representing the underlying state space of the concurrent
system. And an additional structure on this topological space models
time irreversibility.

In the setting of flows introduced in \cite{model3}, the non-constant
execution paths are viewed as objects themselves, not as paths of an
underlying topological space. The topology of the path space models
concurrency and non-constant execution paths can be composed. So a
flow is, by definition, a small category without identity maps
enriched over compactly generated topological spaces.  The main
motivation for introducing this category is the study of the branching
and merging homology theories \cite{exbranch}.  Indeed, they impose
the functoriality of the mapping associating an object with its set of
non-constant execution paths~\footnote{See \cite{model3} Section~20
  for further explanations.} and also the possibility of composing
cubes~\footnote{See the introduction and especially Figure~3 of
  \cite{Coin} for further explanations.}. None of the other
topological models of concurrent systems introduced so far ($d$-space,
local pospace, stream, $d$-space colimit-generated by cubes) has the
first feature since the associated categories contain too many
morphisms. More precisely, they contain morphisms contracting
non-constant execution paths. Of course, it is possible to remove the
``contracting morphisms'' from the categories of $d$-spaces, local
pospaces, streams and $d$-spaces colimit-generated by cubes. But these
models then lose all their interesting properties.  These homology
theories are expected to be important in the study of higher
dimensional bisimulation between concurrent processes by algebraic
invariants.  Indeed, all these notions are related to the structure of
the branching and merging areas of execution paths of a time flow.

It is possible to realize every process algebra \cite{MR1365754} as a
precubical set, and as a flow \cite{ccsprecub} using a realization
functor $|-|_{flow}$. This realization functor is complicated to
handle since its construction requires the use of the cofibrant
replacement functor of the category of flows which is a transfinite
construction of length $2^{\aleph_0}$ (\cite{model3}
Proposition~11.5). The main result of the paper is

\begin{thmN} (Theorem~\ref{small} and Corollary~\ref{small-cor}) There
  exists a \textit{small} realization functor $\gl(-)$ from precubical
  sets to flows which is colimit-preserving. Small means that for
  every $n\geq 0$, there exists a pushout diagram of flows
\[
\xymatrix{
\glob(\mathbf{S}^{n-1}) \fr{} \fd{}& \gl(\de\square[n+1])\fd{}\\
\glob(\mathbf{D}^{n}) \fr{} & \cocartesien\gl(\square[n+1])}
\]
where $\square[n+1]$ is the $(n+1)$-dimensional cube and where
$\de\square[n+1]$ is its boundary.  It is a realization functor in the
sense that there exists a natural transformation $\mu:\gl(-)
\rightarrow |-|_{flow}$ inducing for every precubical set $K$ a
natural S-homotopy equivalence $\mu_K:\gl(K) \simeq |K|_{flow}$ and a
natural transformation $\nu:|-|_{flow} \rightarrow \gl(-)$ inducing
for every precubical set $K$ a natural S-homotopy equivalence
$\nu_K:|K|_{flow} \simeq \gl(K)$ which is an inverse up to S-homotopy
of $\mu_K$.
\end{thmN}

Two applications of this result are given in this paper. Other
applications will be given in future papers.

The lack of a real underlying topological space in the setting of
flows makes some situations very difficult to treat. The first
application is:

\begin{thmN} (Theorem~\ref{unique}) The realization functor
  $|-|_{flow}: \square^{op}\set \rightarrow \dtop$ from the category
  of precubical sets to that of flows defined in \cite{ccsprecub}
  factors up to a natural S-homotopy equivalence as a composite
  \[\square^{op}\set \stackrel{\gl^{top}}\longrightarrow \gltop
  \stackrel{cat}\longrightarrow \dtop\] where $\gltop$ is the category
  of globular complexes and where $cat:\gltop \rightarrow \dtop$ is
  the realization functor from globular complexes to flows defined in
  \cite{model2}. The functor $\gl^{top}$ is unique up to a natural
  S-homotopy equivalence of globular complexes.
\end{thmN}

A notion of underlying homotopy type of flow does exist anyway.  The
underlying state space of a flow exists, and is unique up to homotopy,
not up to homeomorphism \cite{model2}. The fundamental tool to carry
out the construction is also the notion of globular complex.  This
definition enabled us to prove the invariance of the underlying
homotopy type of a flow by refinement of observation in \cite{4eme}.
As a second application of the main result of the paper, or rather, as
an application of the first application, the following theorem
proposes a simplification of the construction of the underlying
homotopy type functor:

\begin{thmN} (Theorem~\ref{app}) Let $K$ be a precubical set. The
  underlying homotopy type of the flow $|K|_{flow}$ associated with
  the precubical set $K$ is naturally isomorphic to the homotopy type
  of the standard cubical complex $|K|_{space}$ associated with $K$:
  i.e.  the functor $| - |_{space}$ is the unique colimit-preserving
  functor from precubical sets to topological spaces associating the
  $n$-cube with the topological $n$-cube $[0,1]^n$ for all $n\geq 0$.
\end{thmN}

This paper can be read as a sequel of the papers \cite{model2} and
\cite{4eme} which study the underlying homotopy type of a flow.
Indeed, several results of \cite{model2} and \cite{4eme} are used in
this work. It can also be read as a continuation of \cite{ccsprecub}
which initializes the study of flows modeling process algebras. This
work is in fact a preparatory work for the study of process algebras
up to homotopy.

\begin{table}
\[
\xymatrix{
\square^{op}\set \ar@{->}[r]^-{\gl^{top}(-)} \ar@/^45pt/[rr]^-{|-|_{flow}} \ar@/^25pt/[rr]^-{\gl(-)}\ar@{->}[d]_-{|-|_{space}}& \gltop  \ar@{->}[dd]^-{\gamma_\gltop}\ar@{->}[ld]^-{|-|}\ar@{->}[r]^-{cat} &
\dtop \ar@{->}[dd]_-{\gamma_\dtop}  \ar@{->} `r[d]
`[dddd] `[ll]_-{\Omega=|-|\circ \overline{cat}^{-1}\circ \gamma_\dtop} [lldd]\\
\top \ar@{->}[d]_-{\gamma_\top} && \\
\ho(\top) & \ar@{->}[l]^-{|-|}\gltop[\mathcal{SH}^{-1}]\iso \gltop/\!\sim_S
\ar@/^20pt/[r]^-{\overline{cat}} \ar@{-}[r]|-{\simeq} &
\ar@/^20pt/[l]^-{\overline{cat}^{-1}}\ho(\dtop)\\
&&\\
&&
}
\]
\begin{tabular}{|c|l|}
  \hline
  $\gamma_\top,\gamma_\gltop,\gamma_\dtop$ & Canonical localization functors\\
  $|-|_{space}$ & Cubical complex associated with a precubical set \\
  $|-|_{flow}$ & Realization functor from precubical sets to  flows\\
  $|-|$ & Underlying topological space functor\\
  $\gl(-)$ & Small realization functor from precubical sets to flows\\
  $\gl^{top}(-)$ & Small realization functor from precubical sets to
  globular complexes\\
  $cat$ and $\overline{cat}^{-1}$ & Equivalence between globular
  complexes and flows\\
  $\Omega$ & Underlying homotopy type functor\\
  \hline
\end{tabular}
\caption{Overview of the functors of the paper}
\end{table}

\subsection{Outline of the paper}

Section~\ref{co} is devoted to the preparatory proofs of some facts
about cocubical objects in a simplicial model category.  The main
Theorem~\ref{homotopy-natural} and Theorem~\ref{homotopy-natural-fini}
are used in the proofs of Theorem~\ref{small}, of
Corollary~\ref{small-cor} and of Theorem~\ref{app}.
Section~\ref{simplicialmodel} constructs the simplicial structure of
the model category of flows. This structure is necessary for the
application of Theorem~\ref{homotopy-natural} and
Theorem~\ref{homotopy-natural-fini} in the proofs of
Theorem~\ref{small} and Corollary~\ref{small-cor}.  This result was
not yet available in a published work.  Section~\ref{section-small}
constructs the small realization functor from precubical sets to
flows. Section~\ref{glob} describes the first application, and
Section~\ref{under} the second application of the main result of the
paper.

\subsection{Prerequisites and notations}

It is required some familiarity with model category techniques
\cite{MR99h:55031} \cite{ref_model2}, with category theory
\cite{MR1712872} \cite{MR96g:18001a} and with simplicial techniques
\cite{MR2001d:55012}. The notation $\simeq$ means \textit{weak
  equivalence} or \textit{equivalence of categories}, the notation
$\iso$ means \textit{isomorphism}, the notation
$\xymatrix@1{\ar@{^{(}->}[r]&}$ means \textit{cofibration} and the
notation $\xymatrix@1{\ar@{->>}[r]&}$ means \textit{fibration}. Let
$\C$ be a cocomplete category.  The class of morphisms of $\C$ that
are transfinite compositions of pushouts of elements of a set of
morphisms $K$ is denoted by $\cell(K)$. An element of $\cell(K)$ is
called a \textit{relative $K$-cell complex}. The category of sets is
denoted by $\set$. The cofibrant replacement functor is denoted by
$(-)^{cof}$.  The function complex of a simplicial model category is
denoted by $\map(-,-)$. The initial object is denoted by
$\varnothing$. The terminal object is denoted by $\mathbf{1}$. In
general, the category of functors from a category $\mathcal{B}$ to a
category $\mathcal{M}$ is denoted by $\mathcal{M}^\mathcal{B}$. Note
that the category $\mathcal{M}^\mathcal{B}$ is locally small if and
only if the category $\mathcal{B}$ is essentially small \cite{small}.
The category of simplicial sets is denoted by $\sis$. $\Delta$ is the
standard category of simplices. $\Delta[n]=\Delta(-,[n])$ is the
standard $n$-simplex.

\section{About cocubical objects in a simplicial model category}
\label{co}

\subsection{Precubical set}

A \textit{precubical set} $K$ consists in a family of sets $(K_n)_{n
  \geq 0}$ and of set maps $\de_i^\alpha:K_n \rightarrow K_{n-1}$ with
$n\geq 1$, $1\leq i \leq n$ and $\alpha\in\{0,1\}$ satisfying the
cubical relations $\de_i^\alpha\de_j^\beta = \de_{j-1}^\beta
\de_i^\alpha$ for any $\alpha,\beta\in \{0,1\}$ and for $i<j$
\cite{Brown_cube}. An element of $K_n$ is called a \textit{$n$-cube}.

A good reference for presheaves is \cite{MR1300636}. A precubical set
can be viewed as a presheaf over a small category denoted by
$\square$~\footnote{All the facts about the small category $\square$
  recalled here are used later in the paper.}  with set of objects
$\{[n],n\in \N\}$, generated by the morphisms $\delta_i^\alpha: [n-1]
\rightarrow [n]$ for $1\leq i\leq n$ and $\alpha \in \{0,1\}$ and
satisfying the cocubical relations $\delta_j^\beta \circ
\delta_i^\alpha = \delta_i^\alpha \circ \delta_{j-1}^\beta $ for $i<j$
and for all $(\alpha,\beta)\in \{0,1\}^2$. With the conventions $[0] =
\{0\}$, $[n] = \{0,1\}^n$ for $n \geq 1$ and $\{0,1\}^0=\{0\}$, the
small category $\square$ is the subcategory of the category of sets
generated by the set maps $\delta_i^\alpha : [n-1] \rightarrow [n]$
for $1\leq i\leq n$ and $\alpha \in \{0,1\}$ defined by
\[\delta_i^\alpha(\epsilon_1, \dots, \epsilon_{n-1}) = (\epsilon_1,
\dots, \epsilon_{i-1}, \alpha, \epsilon_i, \dots, \epsilon_{n-1}).\]
The corresponding category is denoted by $\square^{op}\set$.

Let $\square[n]:=\square(-,[n])$. This defines a functor, also denoted
by $\square$, from $\square$ to $\square^{op}\set$. By Yoneda's lemma,
one has the natural bijection of sets
\[\square^{op}\set(\square[n],K)\iso K_n\] for every precubical set
$K$.  The boundary of $\square[n]$ is the precubical set denoted by
$\de \square[n]$ defined by removing the interior of $\square[n]$:
$(\de \square[n])_k := (\square[n])_k$ for $k<n$ and $(\de
\square[n])_k = \varnothing$ for $k\geq n$.  In particular, one has
$\de \square[0] = \varnothing$.

Let $K$ be a precubical set. Let $K_{\leq n}$ denote the precubical
set obtained from $K$ by keeping the $p$-dimensional cubes of $K$ only
for $p\leq n$. In particular, $K_{\leq 0}=K_0$. Let $\square_n\subset
\square$ be the full subcategory of $\square$ whose set of objects is
$\{[k],k\leq n\}$. A presheaf over $\square_n$ is called a
\textit{$n$-dimensional precubical set}. The category
$\square_n^{op}\set$ will be identified with the full subcategory of
$\square^{op}\set$ of precubical sets $K$ such that the inclusion
$K_{\leq n}\subset K$ is an isomorphism of $\square^{op}\set$.

Let $K$ be a precubical set. The \textit{category $\square \ddownarrow
  K$ of cubes of $K$} is the small category defined by the pullback of
categories
\[
\xymatrix{
\square \ddownarrow K \fr{}\fd{}\cartesien & \square^{op}\set\ddownarrow K
\fd{}\\
\square \fr{} & \square^{op}\set.}
\]
In other terms, an object of $\square \ddownarrow K$ is a morphism
$\square[m]\rightarrow K$ and a morphism of $\square\ddownarrow K$ is
a commutative diagram
\[
\xymatrix{
\square[m] \ar@{->}[rd] \ar@{->}[rr] && \square[n] \ar@{->}[ld]\\
& K.&}
\]

\subsection{The category of all small diagrams over a cocomplete
  category}

Let $\D\mathcal{M}$ be the category of all small diagrams of objects
of a cocomplete category $\mathcal{M}$.  The objects are the functors
$D:\mathcal{B}\rightarrow \mathcal{M}$ where $\mathcal{B}$ is a small
category.  A morphism from a diagram $D:\mathcal{B}\rightarrow
\mathcal{M}$ to a diagram $E:\mathcal{C}\rightarrow \mathcal{M}$ is a
functor $\phi:\mathcal{B}\rightarrow \mathcal{C}$ together with a
natural transformation $\mu:D\rightarrow E\circ \phi$.

\bp \label{fonctoriel} Let $\mathcal{M}$ be a cocomplete category. The
colimit construction $D\mapsto \liminj D$ induces a functor from
$\D\mathcal{M}$ to $\mathcal{M}$. \ep

\bpf A morphism of diagrams $(\phi,\mu):D\rightarrow E$ gives rise to
a morphism $D\rightarrow E\circ \phi$ in $\mathcal{M}^\C$. For every
object $W$ of $\mathcal{M}$, one has the natural set map (where
$\diag_\mathcal{B}$ is the constant diagram functor over $\mathcal{B}$
and where $\diag_\mathcal{C}$ is the constant diagram functor over
$\mathcal{C}$):
\begin{align*}
  & \mathcal{M}(\liminj E,W) & \\
  & \iso \mathcal{M}^\C(E,\diag_\mathcal{C}(W)) & \hbox{ by adjunction}\\
  & \longrightarrow \mathcal{M}^\mathcal{B}(E\circ
  \phi,\diag_\mathcal{C}(W)\circ
  \phi) & \\
  & = \mathcal{M}^\mathcal{B}(E\circ \phi,\diag_\mathcal{B}(W))& \hbox{ since
    $\diag_\mathcal{C}(W)\circ
    \phi=\diag_\mathcal{B}(W)$}\\
    & \longrightarrow \mathcal{M}^\mathcal{B}(D,\diag_\mathcal{B}(W))& \\
    & \iso \mathcal{M}(\liminj D,W) & \hbox{ by adjunction.}
\end{align*}
One obtains a map $\liminj D\longrightarrow \liminj E$ by setting
$W=\liminj E$. So the colimit construction induces a functor from
$\D\mathcal{M}$ to $\mathcal{M}$. \epf

\subsection{Cocubical  objects in a simplicial model category}

Let us consider a simplicial model category $\mathcal{M}$.  A
\textit{cocubical object} (resp. of dimension $n\geq 0$) is a functor
from $\square$ (resp. $\square_n$) to $\mathcal{M}$.

Let $X$ be a cocubical object of $\mathcal{M}$. Let $\widetilde{X}_K$
be the functor from the category of cubes $\square \ddownarrow K$ of a
precubical set $K$ to $\mathcal{M}$ defined on objects by
\[\widetilde{X}_K(\square[n]\rightarrow K)=X([n])\] and on morphisms
by
\[
\widetilde{X}_K\left(
\begin{array}{ccc}
\square[m] & \stackrel{\square[\delta]}\longrightarrow & \square[n] \\
\downarrow && \downarrow\\
K & = & K
\end{array}
\right)=X(\delta).\] The mapping $X\mapsto \widetilde{X}$ induces a
functor from $\mathcal{M}^\square$ to
$\D\mathcal{M}^{\square^{op}\set}$. Let $\widehat{(-)}$ be the composite
functor
\[
\xymatrix{\widehat{(-)}:\mathcal{M}^\square \fr{\widetilde{(-)}} &
  \D\mathcal{M}^{\square^{op}\set} \fr{\liminj}&
  \mathcal{M}^{\square^{op}\set}}\] in which the right-hand functor is
the functor of Proposition~\ref{fonctoriel}. So one has
\[\boxed{\widehat{X}(K)=\liminj_{\square[n]\rightarrow K} X([n]) =
\liminj_{\square K} \widetilde{X}_K}.\]

The category $\square$ has a structure of a direct Reedy category with
the degree function $d([n])=n$ for all $n\geq 0$.  Let us equip the
category $\mathcal{M}^\square$ of cocubical objects of $\mathcal{M}$
with its Reedy model category structure (\cite{ref_model2}
Theorem~15.3.4). The following proposition describes the Reedy
cofibrations and the Reedy fibrations of cocubical objects.

\bp \label{tous} Let $\mathcal{M}$ be a model category. Then: 
\begin{enumerate}
\item The Reedy fibrations of cocubical objects are the objectwise
  fibrations.
\item A cocubical object $X$ of $\mathcal{M}$ is Reedy fibrant if and
  only if for every $n\geq 0$, $X([n])$ is fibrant.
\item A cocubical object $X$ is Reedy cofibrant if and only if for any
  $n\geq 0$, the map $\widehat{X}(\de\square[n]\subset \square[n])$ is
  a cofibration of $\mathcal{M}$. 
\end{enumerate}
\ep 

\bpf A Reedy fibration of $\mathcal{M}^\square$ is by definition a map
of cocubical objects $X\rightarrow Y$ such that for every object $[n]$
of $\square$, the map \[X([n]) \rightarrow M_{[n]}X \p_{M_{[n]}Y}
Y([n])\] is a fibration of $\mathcal{M}$ where $M_{[n]}X$ (resp.
$M_{[n]}Y$) is the matching object of $X$ (resp. of $Y$) at $[n]$.
These matching objects are equal to the terminal object $\mathbf{1}$
of $\mathcal{M}$ since the Reedy category $\square$ is direct. So the
Reedy fibrations are the objectwise fibrations. Hence the first
assertion.

A cocubical object $X$ is therefore Reedy fibrant if and only if for
every object $[n]$ of $\square$, the map $X([n]) \rightarrow
\mathbf{1}$ is a fibration.  Hence the second assertion.

A Reedy cofibration of $\mathcal{M}^\square$ is by definition a map of
cocubical objects $X\rightarrow Y$ such that for every object $[n]$ of
$\square$, the map \[L_{[n]}Y \sqcup_{L_{[n]}X} X([n]) \rightarrow
Y([n])\] is a cofibration of $\mathcal{M}$ where $L_{[n]}X$ (resp.
$L_{[n]}Y$) is the matching object of $X$ (resp. of $Y$) at $[n]$. The
latching category at $\alpha=\square[n]$, usually denoted by
$\de(\square \ddownarrow \alpha)$, is by definition the full
subcategory of the category $\square \ddownarrow \alpha$ of cubes of
$\square[n]$ containing the maps $\square[m]\rightarrow \square[n]$
different from the identity of $\square[n]$. And the latching object
of $X$ at $\alpha$ is by definition \[L_{[n]}X \iso
\liminj_{\de(\square \ddownarrow \alpha)} X(\square[m])\iso
\widehat{X}(\de\square[n]).\] Hence the third assertion.  \epf

\bp \label{equivalence} The composite functor
\[\widehat{(-)}:\mathcal{M}^\square
\stackrel{\widetilde{(-)}}\longrightarrow
\D\mathcal{M}^{\square^{op}\set} \stackrel{\liminj}\longrightarrow
\mathcal{M}^{\square^{op}\set}\] induces an equivalence of categories
\[\mathcal{M}^\square \simeq \mathcal{M}_{\liminj}^{\square^{op}\set}\] 
where $\mathcal{M}_{\liminj}^{\square^{op}\set}$ is the category of
colimit-preserving functors from $\square^{op}\set$ to $\mathcal{M}$.
\ep

\bpf Consider the functor $F:\mathcal{M}_{\liminj}^{\square^{op}\set}
\rightarrow \mathcal{M}^\square$ defined by $F(Z)=Z\circ \square$.
Then for every cocubical object $X$ of $\mathcal{M}$, one has the
natural isomorphisms of $\mathcal{M}$
\[F(\widehat{X})([n])\iso \widehat{X}(\square[n]) \iso
\liminj_{\square[m]\rightarrow \square[n]} X([m]) \iso X([n])\] and,
since $Z$ is colimit-preserving, 
\[ \widehat{F(Z)}(K) \iso \liminj_{\square[n]\rightarrow K}
Z(\square[n]) \iso Z(K).\] 
\epf 

By \cite{ref_model2} Theorem~15.3.4 again, the Reedy model structure
of $\mathcal{M}^\square$ is simplicial with the tensor product and
cotensor product of a cocubical object $X$ by a simplicial set $K$
defined by the composites $X\ot K:= (-\ot K)\circ X$ and $X^ K:= (-)^
K\circ X$.

\bth \label{homotopy-natural} Let $\mathcal{M}$ be a simplicial model
category. Let $I$, $X$ and $Y$ be three cocubical objects of
$\mathcal{M}$. Let $p_X:X\rightarrow I$ and $p_Y:Y\rightarrow I$ be
two objectwise trivial fibrations of cocubical objects of
$\mathcal{M}$.  Assume that for every $n\geq 0$, the maps
$\widehat{X}(\de\square[n]) \rightarrow \widehat{X}(\square[n])$ and
$\widehat{Y}(\de\square[n]) \rightarrow \widehat{Y}(\square[n])$ are
cofibrations of $\mathcal{M}$ and $I([n])$ is fibrant. Then:
\begin{itemize} 
\item There exists a natural transformation from $\widehat{X}$ to
  $\widehat{Y}$ over $\widehat{I}$, i.e. a map $\widehat{X}
  \rightarrow \widehat{Y}$ such that the following diagram commutes:
\[
\xymatrix{
  \widehat{X} \ar@{->}[rd]_-{\widehat{p_X}} \ar@{->}[rr] && \widehat{Y} \ar@{->}[ld]^-{\widehat{p_Y}}\\
  & \widehat{I}.&}
\]
\item Take two natural transformations $\widehat{\mu}:\widehat{X}
  \rightarrow \widehat{Y}$ and $\widehat{\nu}:\widehat{X} \rightarrow
  \widehat{Y}$ over $\widehat{I}$. Then there exists a simplicial
  homotopy between $\widehat{\mu}(K)$ and $\widehat{\nu}(K)$ which is
  natural with respect to $K$.
\item For any natural transformation $\widehat{\mu}:\widehat{X}
  \rightarrow \widehat{Y}$ over $\widehat{I}$ and any natural
  transformation $\widehat{\nu}:\widehat{Y} \rightarrow \widehat{X}$
  over $\widehat{I}$, the map $\widehat{\mu}(K)\circ \widehat{\nu}(K)$
  is naturally simplicially homotopy equivalent to
  $\id_{\widehat{Y}(K)}$ and the map $\widehat{\nu}(K)\circ
  \widehat{\mu}(K)$ is naturally simplicially homotopy equivalent to
  $\id_{\widehat{X}(K)}$, natural meaning natural with respect to $K$.
\end{itemize}
\eth 

\bpf The cocubical object $X$ is Reedy cofibrant by
Proposition~\ref{tous}.  The map $p_Y$ is a trivial Reedy fibration by
Proposition~\ref{tous} as well. Let $A\rightarrow B$ be a cofibration
of simplicial sets. Consider a commutative diagram of simplicial sets:
\[
\xymatrix{
A \fr{} \fd{} & \map(X,Y) \fd{}\\
B \fr{} \ar@{-->}[ru]^-{k}& \map(X,I)}
\] 
where the map $(p_Y)_*:\map(X,Y) \rightarrow \map(X,I)$ is induced by
the composition by $p_Y$. By adjunction, the lift $k$ exists if and
only if the lift $k'$ exists in the commutative diagram of cocubical
objects of $\mathcal{M}$
\[
\xymatrix{
X\ot A \fr{} \fd{} & Y \fd{}\\
X\ot B \fr{} \ar@{-->}[ru]^-{k'}& I.}
\]
The map $X\ot A \rightarrow X\ot B$ is the pushout product of the
Reedy cofibration $\varnothing \rightarrow X$ by the cofibration of
simplicial sets $A\rightarrow B$, and therefore is a Reedy
cofibration.  Hence the existence of $k'$ and $k$. So the simplicial
map $(p_Y)_*:\map(X,Y) \rightarrow \map(X,I)$ is a trivial fibration
of simplicial sets. Let $F$ be the fibre over $p_X$ of the simplicial
map $(p_Y)_*:\map(X,Y) \rightarrow \map(X,I)$ defined by the pullback
diagram of simplicial sets:
\[
\xymatrix{
  F\cartesien \fr{} \ar@{->>}[d]^-{\simeq} & \map(X,Y) \ar@{->>}[d]^-{\simeq} \\
  \Delta[0] \fr{p_X} & \map(X,I).}
\]
Since the pullback of a trivial fibration is a trivial fibration, the
map $F\rightarrow \Delta[0]$ is a trivial fibration. The lift $\ell$
in the commutative diagram of simplicial sets
\[
\xymatrix{
\varnothing \fr{} \ar@{^{(}->}[d] &  F  \ar@{->>}[d]^-{\simeq} \\
\Delta[0] \ar@{=}[r]\ar@{-->}[ru]^-{\ell} & \Delta[0]  }
\]
gives $\ell(0)\in F_0 \subset \map(X,Y)_0=\mathcal{M}^\square(X,Y)$.
By definition, $\ell(0):X\rightarrow Y$ is a morphism of cocubical
objects over $I$.  Hence a natural transformation
$\widehat{\ell(0)}:\widehat{X} \rightarrow \widehat{Y}$ over
$\widehat{I}$ and the first assertion.

Take two natural transformations $\widehat{\mu}$ and $\widehat{\nu}$
from $\widehat{X}$ to $\widehat{Y}$ over $\widehat{I}$. By
Proposition~\ref{equivalence}, one can suppose that they come from two
morphisms of cocubical objects ${\mu}$ and ${\nu}$ from ${X}$ to ${Y}$
over $I$. One obtains the commutative diagram of simplicial sets:
\[
\xymatrix{
\Delta[0] \sqcup \Delta[0]\ar@{^{(}->}[d] \fr{(\mu,\nu)} & F \ar@{->>}[d]^-{\simeq}\\
\Delta[1]  \fr{} \ar@{-->}[ru]& \Delta[0]}
\] 
Thus, there exists a simplicial path $\Delta[1] \rightarrow F\subset
\map(X,Y)$ between ${\mu}$ and ${\nu}$, i.e. by adjunction a morphism
of cocubical objects $H:X\ot \Delta[1] \rightarrow Y$ such that the
two natural transformations $X\rightrightarrows X\ot \Delta[1]
\rightarrow Y$ are $\mu$ and $\nu$, i.e. $H$ is a simplicial homotopy
between $\mu$ and $\nu$.  One obtains a simplicial homotopy
$\widetilde{H}_K\in \map(\widetilde{X}_K,\widetilde{Y}_K)_1$ between
$\widetilde{\mu}_K$ and $\widetilde{\nu}_K$. Since there is an
isomorphism \[\liminj (\widetilde{X}_K\ot \Delta[1]) \iso
(\liminj\widetilde{X}_K)\ot \Delta[1]\] because the functor $-\ot
\Delta[1]$ is colimit-preserving, one obtains a simplicial homotopy
$\widehat{H}(K)\in \map(\widehat{X}(K),\widehat{Y}(K))_1$ between
$\widehat{\mu}(K)$ and $\widehat{\nu}(K)$ which is natural with
respect to $K$. Hence the second assertion.

The third assertion is a consequence of the second assertion by
noticing that $\id_{\widehat{X}}$ is a natural transformation from
$\widehat{X}$ to itself over $\widehat{I}$ and that
$\id_{\widehat{Y}}$ is a natural transformation from $\widehat{Y}$ to
itself over $\widehat{I}$.  \epf

\bth \label{homotopy-natural-fini} Let $\mathcal{M}$ be a
  simplicial model category. Let $n\geq 0$. Let $I$, $X$ and $Y$ be
  three cocubical objects of $\mathcal{M}$ of dimension $n$. Let
  $p_X:X\rightarrow I$ and $p_Y:Y\rightarrow I$ be two objectwise
  trivial fibrations of cocubical objects of $\mathcal{M}$.  Assume
  that for every $0 \leq p\leq n$, the maps
  $\widehat{X}(\de\square[p]) \rightarrow \widehat{X}(\square[p])$ and
  $\widehat{Y}(\de\square[p]) \rightarrow \widehat{Y}(\square[p])$ are
  cofibrations of $\mathcal{M}$ and $I([p])$ is fibrant. Then:
\begin{itemize} 
\item There exists a natural transformation from $\widehat{X}$ to
  $\widehat{Y}$ over $\widehat{I}$.
\item Take two natural transformations $\widehat{\mu}:\widehat{X}
  \rightarrow \widehat{Y}$ and $\widehat{\nu}:\widehat{X} \rightarrow
  \widehat{Y}$ over $\widehat{I}$. Then there exists a simplicial
  homotopy between $\widehat{\mu}(K)$ and $\widehat{\nu}(K)$ which is
  natural with respect to $K$.
\item For any natural transformation $\widehat{\mu}:\widehat{X}
  \rightarrow \widehat{Y}$ over $\widehat{I}$ and any natural
  transformation $\widehat{\nu}:\widehat{Y} \rightarrow \widehat{X}$
  over $\widehat{I}$, the map $\widehat{\mu}(K)\circ \widehat{\nu}(K)$
  is naturally simplicially homotopy equivalent to
  $\id_{\widehat{Y}(K)}$ and the map $\widehat{\nu}(K)\circ
  \widehat{\mu}(K)$ is naturally simplicially homotopy equivalent to
  $\id_{\widehat{X}(K)}$, natural meaning natural with respect to the
  precubical set $K$ of dimension $n$ .
\end{itemize}
\eth

\bpf Use the Reedy structure of $\square_n$ and the Reedy model
structure of $\mathcal{M}^{\square_n}$ in the proof of
Theorem~\ref{homotopy-natural}. \epf

\section{The weak S-homotopy model category of flows is simplicial}
\label{simplicialmodel}

The goal of this section is to prove that the weak S-homotopy model
category of $\dtop$ is simplicial (Theorem~\ref{simpl}).

\subsection{Topological space}

All topological spaces are compactly generated, i.e. weak Hausdorff
$k$-spaces. Further details about these topological spaces are
available in \cite{MR90k:54001} \cite{MR2000h:55002}, the appendix of
\cite{Ref_wH} and also the preliminaries of \cite{model3}. All compact
spaces are Hausdorff.  The category of \textit{compactly generated
  topological space}s together with the continuous maps is denoted by
$\top$.  The category $\top$ is equipped with the usual model
structure having the weak homotopy equivalences as weak equivalences
and having the Serre fibrations as fibrations. This model structure is
simplicial and any topological space is fibrant. The homotopy
category, i.e. the localization of $\top$ by the weak homotopy
equivalences, is denoted by $\ho(\top)$. The functor $\gamma_\top:
\top \rightarrow \ho(\top)$ is the canonical functor which is the
identity on objects. The set of continuous maps $\top(X,Y)$ from $X$
to $Y$ equipped with the Kelleyfication of the compact-open topology
is denoted by $\ttop(X,Y)$. The latter topological space is the
internal hom of the cartesian closed category $\top$. So one has the
natural homeomorphism $\ttop(X\p Y,Z) \iso \ttop(X,\ttop(Y,Z))$.
Moreover the covariant functor $\ttop(X,-):\top\rightarrow \top$ and
the contravariant functor $\ttop(-,X):\top^{op}\rightarrow \top$
preserve limits \cite{MR651714}.

\subsection{Flow}

A \textit{flow} $X$ is a small category without identity maps enriched
over the category of compactly generated topological spaces. The set
of objects is denoted by $X^0$. The space of morphisms from $\alpha$
to $\beta$ is denoted by $\P_{\alpha,\beta} X$~\footnote{Sometimes, an
  object of a flow is called a state and a morphism a (non-constant)
  execution path.}. A morphism of flows $f:X \rightarrow Y$ is a set
map $f^0:X^0 \rightarrow Y^0$ together with a continuous map $\P f:\P
X \rightarrow \P Y$ preserving the structure. The corresponding
category is denoted by $\dtop$.  If for all $\alpha \in X^0$, the
space $\P_{\alpha,\alpha}X$ is empty, then $X$ is called a
\textit{loopless flow}. The composition law of a flow is denoted by
$*$.

Any poset $P$, and in particular any set, can be viewed as a loopless
flow with a morphism from $\alpha$ to $\beta$ if and only if
$\alpha<\beta$. This yields a functor from the category of posets with
strictly increasing maps to that of flows. 

Let $Z$ be a topological space. The flow $\glob(Z)$ is defined by
\begin{itemize}
\item $\glob(Z)^0=\{\widehat{0},\widehat{1}\}$, 
\item $\P \glob(Z)= \P_{\widehat{0},\widehat{1}} \glob(Z) = Z$,
\item $s=\widehat{0}$, $t=\widehat{1}$ and a trivial composition law.
\end{itemize}
It is called the \textit{globe} of the space $Z$.

The weak S-homotopy model structure of $\dtop$ is characterized as
follows \cite{model3}:
\begin{itemize}
\item The weak equivalences are the \textit{weak S-homotopy
    equivalences}, i.e. the morphisms of flows $f:X\longrightarrow Y$
  such that $f^0:X^0\longrightarrow Y^0$ is a bijection of sets and
  such that $\P f:\P X\longrightarrow \P Y$ is a weak homotopy
  equivalence.
\item The fibrations are the morphisms of flows
$f:X\longrightarrow Y$ such that $\P f:\P X\longrightarrow \P Y$ is a
Serre fibration. 
\end{itemize}
Note that any flow is fibrant.  The homotopy category is denoted by
$\ho(\dtop)$. This model structure is cofibrantly generated. The set
of generating cofibrations is the set \[I^{gl}_{+}=I^{gl}\cup
\{R:\{0,1\}\longrightarrow \{0\},C:\varnothing\longrightarrow
\{0\}\}\] with $I^{gl}=\{\glob(\mathbf{S}^{n-1})\subset
\glob(\mathbf{D}^{n}), n\geq 0\}$ where $\mathbf{D}^{n}$ is the
$n$-dimensional disk and $\mathbf{S}^{n-1}$ the $(n-1)$-dimensional
sphere. By convention, the $(-1)$-dimensional sphere is the empty
space. The set of generating trivial cofibrations is
\[J^{gl}=\{\glob(\mathbf{D}^{n}\p\{0\})\subset \glob(\mathbf{D}^{n}\p
[0,1]), n\geq 0\}.\]

Let $X$ and $U$ be two flows. Let $\tdtop(X,U)$ be the set
$\dtop(X,U)$ equipped with the Kelleyfication of the compact-open
topology. Let $f,g:X\rightrightarrows U$ be two morphisms of flows.
Then a \textit{S-homotopy} is a continuous map $H:[0,1] \rightarrow
\tdtop(X,U)$ with $H_0=f$ and $H_1=g$. This situation is denoted by
$f\sim_S g$.  The S-homotopy relation defines a congruence on the
category $\dtop$.  If there exists a map $f':U\rightarrow X$ with
$f\circ f'\sim_S \id_U$ and $f'\circ f\sim_S \id_X$, then $f$ is
called a \textit{S-homotopy equivalence}.

\subsection{Simplicial structure of the model category of flows}

Let $K$ be a non-empty connected simplicial set. Let $X$ be a flow.
Then one has the isomorphism of topological spaces
\[\ttop(|K|,\P X) \iso \bigsqcup_{(\alpha,\beta)\in X^0\p
  X^0}\ttop(|K|,\P_{\alpha,\beta}X)\] where the topological space
$|K|$ is the geometric realization of the simplicial set $K$. Note the
latter isomorphism is false if $K$ is empty or not connected.  The
associative composition law $*:\P X\p_{X^0} \P X \rightarrow \P X$
gives rise to a continuous map \[*: \ttop(|K|,\P X)\p_{X^0}
\ttop(|K|,\P X) \iso \ttop(|K|,\P X\p_{X^0} \P X) \rightarrow
\ttop(|K|,\P X).\] The homeomorphism $\ttop(|K|,\P X)\p_{X^0}
\ttop(|K|,\P X) \iso \ttop(|K|,\P X\p_{X^0} \P X)$ holds since the
functor $\ttop(|K|,-)$ is limit-preserving. Hence the following
definition:

\bd Let $K$ be a non-empty connected simplicial set. Let $X$ be an
object of $\dtop$. Let $X^K$ be the flow defined by
\begin{itemize}
\item $(X^K)^0=X^0$
\item $\P_{\alpha,\beta} (X^K)=\ttop(|K|,\P_{\alpha,\beta} X)$ for all
  $(\alpha,\beta)\in X^0\p X^0$
\item the above composition law.
\end{itemize}
\ed 

Several theorems of \cite{model3} are going to be used. Therefore, it
is helpful for the reader to give the following correspondence between
the notations of this paper (left column) and of \cite{model3} (right
column) for a flow $X$ and a non-empty connected simplicial set $K$:
\begin{center}
\begin{tabular}{|c|c|}
\hline 
Notations of this paper & Notations of \cite{model3} \\
\hline 
$X \ot K$ & $|K| \boxtimes X$ \\
\hline 
$X^K$ & $\{|K|,X\}_S$\\
\hline 
\end{tabular}
\end{center}

\bp \label{adj1moins} Let $K$ be a non-empty connected simplicial set.
The mapping $X \mapsto X^K$ gives rise to an undofunctor of $\dtop$.
The functor $(-)^K$ is a right adjoint.  \ep

\bpf Consequence of \cite{model3} Theorem~7.8.  \epf

\begin{nota} Let $K$ be a non-empty connected simplicial set. Let us
  denote by $-\ot K$ the left adjoint of $(-)^K$.
\end{nota}

\bd Let $K$ be a non-empty simplicial set. Let $(K_i)_{i\in I}$ be its
set of non-empty connected components. Let 
\begin{itemize} 
\item $X \ot K := \bigsqcup_{i\in I} X \ot K_i$
\item $X^K := \prod_{i\in I} X^{K_i}$. 
\end{itemize}
And let $X \ot \varnothing = \varnothing$ and $X ^\varnothing =
\mathbf{1}$.  \ed

\bp \label{adj2moins} Let $K$ be a simplicial set. The pair of
functors $(- \ot K):\dtop \leftrightarrows \dtop:(-)^K$ is a
categorical adjunction.  \ep

\bpf If $K=\varnothing$, then one has the isomorphisms 
\[\dtop(X\ot \varnothing,Y)\iso \dtop(\varnothing,Y) \iso
\mathbf{1} \iso \dtop(X,Y^\varnothing).\] Now let $K$ be a non-empty
simplicial set. Let $(K_i)_{i\in I}$ be its set of non-empty connected
components. Then one has 
\begin{align*}
& \dtop(X\ot K,Y) & \\
& \iso \dtop(\bigsqcup_{i\in I} (X\ot K_i),Y) & \hbox{ by definition of $X\ot K$}\\
& \iso \prod_{i\in I} \dtop(X\ot K_i,Y) & \\
& \iso \prod_{i\in I} \dtop(X,Y^{K_i}) & \hbox{ by Proposition~\ref{adj1moins}}\\
& \iso  \dtop(X, \prod_{i\in I}Y^{K_i}) & \\
& \iso \dtop(X,Y^K) &   \hbox{ by definition of $X^K$.}
\end{align*}
Hence the adjunction.  \epf

\bd Let $X$ and $Y$ be two objects of $\dtop$. Let
\[\map(X,Y) := \dtop(X\ot \Delta[*],Y).\] 
It is called the {\rm function complex} from $X$ to $Y$ \ed

\bp \label{singmap} Let $X$ and $Y$ be two flows. Then one has the natural 
isomorphism of simplicial sets 
\[\map(X,Y) \iso \sing (\tdtop(X,Y))\] 
where $\sing$ is the singular nerve functor.  \ep

\bpf Let $n\geq 0$. Since the topological space $|\Delta[n]$ is
non-empty and connected, one has $\sing (\tdtop(X,Y))_n =
\top(|\Delta[n]|,\tdtop(X,Y)) \iso \dtop(X \ot \Delta[n],Y)$ by
\cite{model3} Theorem~7.9.  \epf

\bp \label{limcolimmoins} Let $\mathcal{B}$ be a small category. Let
$X:\mathcal{B}\rightarrow \dtop$ be a functor. Then one has the
natural isomorphisms of simplicial sets \beas
&& \map(\liminj X_b,Y) \iso \limproj \map(X_b,Y) \\
&& \map(Y,\limproj X_b) \iso \limproj \map(Y,X_b) \eeas for any flow
$Y$ of $\dtop$.  \ep

\bpf Limits and colimits are calculated pointwise in the category of
simplicial sets. Since for every $n\geq 0$, the functor $-\ot
\Delta[n]$ is a left adjoint by Proposition~\ref{adj2moins}, one
obtains the natural bijections 
\begin{align*}
& \map(\liminj X_b,Y)_n & \\
&\iso \dtop((\liminj X_b)\ot \Delta[n],Y) &\hbox{ by definition of $\map$}\\
&\iso \dtop(\liminj (X_b\ot \Delta[n]),Y) & \hbox{ since $-\ot
  \Delta[n]$ is a left adjoint}\\
&\iso \limproj \dtop(X_b\ot \Delta[n],Y) & \\
&\iso \limproj \map(X_b,Y)_n & \hbox{ by definition of $\map$.}
\end{align*}
and 
\begin{align*}
&\map(Y,\limproj X_b)_n & \\
&\iso \dtop(Y\ot \Delta[n],\limproj X_b) &\hbox{ by definition of
  $\map$}\\
&\iso \limproj \dtop(Y\ot \Delta[n],X_b)& \\
& \iso \limproj \map(Y,X_b)_n & \hbox{ by definition of $\map$.}
\end{align*}
\epf

\bp \label{prodprodmoins} Let $K$ and $L$ be two simplicial sets. Then
one has a natural isomorphism of flows $X^{K\p L} = (X^K)^L$ for every
flow $X$ of $\dtop$. \ep

\bpf If $K$ or $L$ is empty, then $X^{K\p L} = (X^K)^L=\mathbf{1}$ by
definition.  Suppose now that $K$ and $L$ are both non-empty and
connected.  The flows $X^{K\p L}$ and $(X^K)^L$ have same set of
states $X^0$.  And $\P(X^{K\p L}) \iso \ttop(|K\p L|,\P X)$ and
$\P((X^K)^L)=\ttop(|L|,\ttop(|K|,\P X))$. Hence the conclusion in this
case since $\top$ is cartesian closed and since there is a
homeomorphism $|K\p L|\iso [K|\p |L|$.  We treat now the general case
where $K$ and $L$ are both non-empty.  Let $(K_i)_{i\in I}$ and
$(L_j)_{j\in J}$ be the non-empty connected components of $K$ and $L$
respectively. Then $K\p L=\prod_{i\in I} \prod_{j\in J} K_i\p L_j$
and:
\begin{align*}
& X^{K\p L} & \\
& \iso \prod_{i\in I} 
\prod_{j\in J} X^{K_i\p L_j} & \hbox{ since the $K_i\p L_j$'s are
  non-empty and connected} \\
& \iso \prod_{i\in I} 
\prod_{j\in J} (X^{K_i})^{L_j} & \\
& \iso \prod_{j\in J} (\prod_{i\in I} X^{K_i})^{L_j} & \hbox{ since the functors $(-)^{L_j}$ are right adjoints}\\
& \iso \prod_{j\in J} (X^K)^{L_j} & \hbox{ by definition of $X^K$}\\
& \iso (X^K)^L & \hbox{ by definition of $(-)^L$.}
\end{align*}
\epf 

\bp \label{ototmoins} Let $K$ and $L$ be two simplicial sets. Let $X$
be a flow. Then one has a natural isomorphism of flows $(X\ot K) \ot
L\iso X\ot (K\p L)$. \ep

\bpf Let $Y$ be another flow. Then one has 
\begin{align*}
  & \dtop((X\ot K) \ot L,Y) & \\
  & \iso \dtop(X\ot K,Y^L) & \hbox{ by Proposition~\ref{adj2moins}}\\
  & \iso \dtop(X,(Y^L)^K) & \hbox{ by Proposition~\ref{adj2moins}} \\
  & \iso \dtop(X,Y^{L\p K}) & \hbox{ by Proposition~\ref{prodprodmoins}}\\
  & \iso \dtop(X\ot (K\p L), Y) & \hbox{ by
    Proposition~\ref{adj2moins}.}
\end{align*}
Hence the result using Yoneda's lemma. 
\epf 

\bp\label{closeness}
Let $K$ be a simplicial set. Let $X$ and $Y$ be two flows. Then one
has a natural isomorphism of simplicial sets
\[\map(X\otimes K,Y)\iso \map(K,\map(X,Y)).\]
\ep

\bpf If $K=\varnothing$, then one has to compare $\map(X\otimes
\varnothing,Y) \iso \map(\varnothing,Y) \iso \mathbf{1}$ by
Proposition~\ref{limcolimmoins} and $\map(\varnothing,\map(X,Y)) \iso
\mathbf{1}$. So one can suppose the simplicial set $K$ non-empty. By
construction of the functor $-\ot K$ and by
Proposition~\ref{limcolimmoins}, one can suppose that $K$ is connected
as well. Let $n\geq 0$.  Thus,
\begin{align*}
  & \map(X\otimes K,Y)_n & \\
  & \iso \dtop(X \otimes (K \p \Delta[n]),Y) & \hbox{ by definition of
    $\map$ and by Proposition~\ref{ototmoins}}\\
  & \iso \top(|K\p \Delta[n]|,\tdtop(X,Y)) & \hbox{ by \cite{model3} Theorem~7.9}\\
  & \iso \sis(K\p \Delta[n],\map(X,Y)) & \hbox{by adjunction and by Proposition~\ref{singmap}}\\
  &\iso \map(K,\map(X,Y))_n& \hbox{ by definition of $\map$ in
    $\sis$.}
\end{align*}
\epf

\bp \label{adj3moins} One has the natural isomorphism of simplicial sets
\[\map(X \ot K,Y) \iso \map (X,Y^K)\]  
for every simplicial set $K$ and every flow $X,Y$ of $\dtop$.  \ep

\bpf
One has for any $n\geq 0$
\begin{align*}
  & \map(X\ot K,Y)_n & \\
  & \iso \dtop(X \ot K \ot \Delta[n], Y) & \hbox{ by definition of
    $\map$}\\
  & \iso \dtop(X \ot \Delta[n],Y^K) & \hbox{ by Proposition~\ref{ototmoins} and Proposition~\ref{adj2moins}}\\
  & \iso \map(X,Y^K)_n. & \hbox{ by definition of $\map$.}
\end{align*}
\epf

\begin{lem} \label{l1}
Let $f:X\longrightarrow Y$ be a morphism of flows. Then the following
conditions are equivalent:
\begin{enumerate}
\item $f$ is a fibration of flows, that is the continuous map $\P f:\P X\longrightarrow \P Y$
is a fibration of topological spaces.
\item for any $(\alpha,\beta)\in X^0\p X^0$, the continuous map
$\P f:\P_{\alpha,\beta}X\longrightarrow \P_{f(\alpha),f(\beta)}Y$
is a fibration of topological spaces.
\end{enumerate}
\end{lem}

\bpf A continuous map is a fibration if and only if it satisfies the
right lifting property with respect to the inclusion maps
$\mathbf{D}^n\longrightarrow \mathbf{D}^n\p [0,1]$. The result comes
from the connectedness of both $\mathbf{D}^n$ and $\mathbf{D}^n\p
[0,1]$.  \epf

\bp \label{fin} Let $i:A\longrightarrow B$ be a cofibration of flows.
Let $p:X\longrightarrow Y$ be a fibration of flows. Then the morphism
of simplicial sets
\[Q(i,p):\map(B,X) \longrightarrow \map(A,X)\p_{\map(A,Y)}\map(B,Y)\]
is a fibration of simplicial sets. Moreover if either $i$ or $p$ is
trivial, then the fibration $Q(i,p)$ is trivial as well. \ep

\bpf By \cite{MR2001d:55012} Proposition~II.3.13 p95, it suffices to
prove that the morphism of flows
\[X^{\Delta[n]}\longrightarrow X^{\de\Delta[n]} \p_{Y^{\de\Delta[n]}}
Y^{\Delta[n]}\] is a fibration (resp. trivial fibration) for any
$n\geq 0$ as soon as $X\longrightarrow Y$ is a fibration (resp.
trivial fibration) of flows. The case of a fibration is the only one
treated since the other case is similar.

For $n=0$, one has to check that
\[X^{\Delta[0]}\longrightarrow \mathbf{1} \p_{\mathbf{1}}
Y^{\Delta[0]}\iso Y^{\Delta[0]}\] is a fibration. Since $\Delta[0]$ is
connected, $X^{\Delta[0]}=X$ and $Y^{\Delta[0]}=Y$.  So there is
nothing to check for $n=0$.

For $n=1$, one has to check that
\[X^{\Delta[1]}\longrightarrow X^{\de\Delta[1]} \p_{Y^{\de\Delta[1]}}
Y^{\Delta[1]}\] is a fibration. Since $\de\Delta[1]$ is the discrete
two-point simplicial set, then $X^{\de\Delta[1]}=X\p X$. Since
$\Delta[1]$ is connected, one has to check that for any
$(\alpha,\beta)\in X^0\p X^0$, the continuous map
\[\ttop(|\Delta[1]|,\P_{\alpha,\beta}X)\longrightarrow
(\P_{\alpha,\beta}X\p \P_{\alpha,\beta}X)\p_{(\P_{\alpha,\beta}Y\p \P_{\alpha,\beta}Y)}
\ttop(|\Delta[1]|,\P_{\alpha,\beta}Y)\]
is a fibration of topological spaces. Using the homeomorphisms
\[\P_{\alpha,\beta}X\p \P_{\alpha,\beta}X
\iso \ttop(\{-1,1\},\P_{\alpha,\beta}X)\] and
\[\P_{\alpha,\beta}Y\p \P_{\alpha,\beta}Y
\iso \ttop(\{-1,1\},\P_{\alpha,\beta}Y),\] one has to check that the
continuous map
\[\ttop(|\Delta[1]|,\P_{\alpha,\beta}X)\longrightarrow
\ttop(\{-1,1\},\P_{\alpha,\beta}X)\p_{\ttop(\{-1,1\},\P_{\alpha,\beta}Y)}
\ttop(|\Delta[1]|,\P_{\alpha,\beta}Y)\] is a fibration of topological
spaces. So one has to prove that for any commutative square 
\[
\xymatrix{
\mathbf{D}^n\p\{0\} \fr{}\fd{} &
\ttop(|\Delta[1]|,\P_{\alpha,\beta}X)\fd{}\\
\mathbf{D}^n\p[0,1] \fr{} \ar@{-->}[ru]^-{k}& \ttop(\{-1,1\},\P_{\alpha,\beta}X)\p_{\ttop(\{-1,1\},\P_{\alpha,\beta}Y)}
\ttop(|\Delta[1]|,\P_{\alpha,\beta}Y),}
\] 
the lift $k$ exists. By adjunction, it suffice to prove that the lift
$k'$ of the commutative square
\[
\xymatrix{ (\mathbf{D}^n\p |\Delta[1]|)\cup (\mathbf{D}^n\p[0,1]\p
  \{-1,1\}) \fr{}\fd{} &
  \P_{\alpha,\beta}X\fd{}\\
  \mathbf{D}^n\p[0,1]\p |\Delta[1]| \fr{} \ar@{-->}[ru]_-{k'}& \P_{\alpha,\beta}Y}
\] 
exists. The inclusion $\{-1,1\}\subset |\Delta[1]|$ is a cofibration.
So the left-hand map is the pushout product of a cofibration with a
trivial cofibration. By Lemma~\ref{l1}, the continuous map
$\P_{\alpha,\beta}X\longrightarrow \P_{\alpha,\beta}Y$ is a fibration
of topological spaces. The case $n=1$ is therefore solved.

Consider now the case $n\geq 2$. Then both $\Delta[n]$ and
$\de\Delta[n]$ are connected. Therefore one has to check that for any
$(\alpha,\beta)\in X^0\p X^0$, the continuous map
\[\ttop(|\Delta[n]|,\P_{\alpha,\beta}X)\longrightarrow \ttop(|\de\Delta[n]|,\P_{\alpha,\beta}X)\p_{\ttop(|\de\Delta[n]|,\P_{\alpha,\beta}Y)}\ttop(|\Delta[n]|,\P_{\alpha,\beta}Y)\]
is a fibration. This holds for the same reason as above because 1) the
category of topological spaces is cartesian closed, 2) the inclusion
$|\de\Delta[n]|\longrightarrow |\Delta[n]|$ is a cofibration and 3)
the mapping $\P_{\alpha,\beta}X\longrightarrow \P_{\alpha,\beta}Y$ is
a fibration by Lemma~\ref{l1}.  \epf

\bth\label{simpl} The model category $\dtop$ together with the
functors $-\ot K$, $(-)^K$ and $\map(-,-)$ assembles to a simplicial
model category.  \eth

\bpf By definition, the set $\map(X,Y)_0$ is the set $\dtop(X,Y)$ of
morphisms from the flow $X$ to the flow $Y$. Thus, the identity of
$\id_X$ yields for every flow $X$ a simplicial map $\Delta[0]
\rightarrow \map(X,X)$. The theorem is then a consequence of
Proposition~\ref{closeness}, Proposition~\ref{adj3moins} and
Proposition~\ref{fin}.  \epf

Let us conclude this section by an important fact:

\bp \label{connexion} Let $X$ and $Y$ be two flows.  Let
$f,g:X\rightrightarrows Y$ be two morphisms of flows. There exists a
simplicial homotopy $H:X\ot \Delta[1] \rightarrow Y$ between $f$ and
$g$ if and only if there exists a S-homotopy $\overline{H}:[0,1]
\rightarrow \tdtop(X,Y)$ between $f$ and $g$. \ep

\bpf Since $\Delta[1]$ is non-empty and connected, one has the
equality $X\ot \Delta[1] = |\Delta[1]| \boxtimes X$ by definition of
the tensor product. And one has the bijection $\dtop( |\Delta[1]|
\boxtimes X,Y)\iso \top([0,1],\tdtop(X,Y)$ by \cite{model3}
Theorem~7.9.  \epf

\section{Comparing realization functors from precubical sets to flows}
\label{section-small}

\subsection{Realizing a precubical set as a flow}

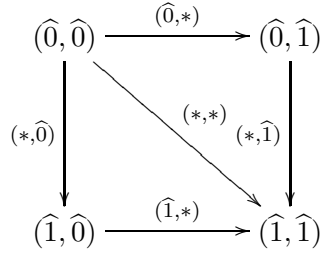
\begin{figure}
\[
\xymatrix{
(\widehat{0},\widehat{0}) \ar@{->}[rr]^-{(\widehat{0},*)}\ar@{->}[dd]_-{(*,\widehat{0})}\ar@{->}[ddrr]^-{(*,*)} && (\widehat{0},\widehat{1})\ar@{->}[dd]_-{(*,\widehat{1})}\\
&& \\
(\widehat{1},\widehat{0}) \ar@{->}[rr]^-{(\widehat{1},*)}&& (\widehat{1},\widehat{1})}
\]
\caption{The flow $\{\widehat{0}<\widehat{1}\}^2$}
\label{cube2}
\end{figure}

A state of the flow associated with the poset
$\{\widehat{0}<\widehat{1}\}^n$ (i.e. the product of $n$ copies of
$\{\widehat{0}<\widehat{1}\}$) is denoted by a $n$-uple of elements of
$\{\widehat{0},\widehat{1}\}$. By convention,
$\{\widehat{0}<\widehat{1}\}^0=\{0\}$. The unique morphism/execution
path from $(x_1,\dots,x_n)$ to $(y_1,\dots,y_n)$ is denoted by a
$n$-uple $(z_1,\dots,z_n)$ of $\{\widehat{0},\widehat{1},*\}$ with
$z_i=x_i$ if $x_i=y_i$ and $z_i=*$ if $x_i<y_i$. For example in the
flow $\{\widehat{0}<\widehat{1}\}^2$ (cf. Figure~\ref{cube2}), one has
the algebraic relation $(*,*) = (\widehat{0},*)*(*,\widehat{1}) =
(*,\widehat{0}) * (\widehat{1},*)$.

Let $\square \rightarrow \dtop$ be the functor defined on objects by
the mapping $[n]\mapsto (\{\widehat{0}<\widehat{1}\}^n)^{cof}$ and on
morphisms by the mapping $ \delta_i^\alpha \mapsto \lp (\epsilon_1,
\dots, \epsilon_{n-1}) \mapsto (\epsilon_1, \dots, \epsilon_{i-1},
\alpha, \epsilon_i, \dots, \epsilon_{n-1})\rp^{cof}$ where the
$\epsilon_i$'s are elements of $\{\widehat{0},\widehat{1},*\}$.

The functor $| - |_{flow}: \square^{op}\set \rightarrow \dtop$ is then
defined by \cite{ccsprecub} \[\boxed{|K|_{flow}:= \liminj_{\square[n] \rightarrow K}
(\{\widehat{0}<\widehat{1}\}^n)^{cof}}.\] It is a left adjoint. So it
commutes with all small colimits. 

The functor $X\mapsto X^0$ from flows to sets is a left adjoint since
there is a natural bijection $\set(X^0,S)\iso \dtop(X,\widehat{S})$
where $\widehat{S}$ is the flow defined by $\widehat{S}^0=S$ and
$\P_{\alpha,\beta}\widehat{S}=\{(\alpha,\beta)\}$ with the composition
law $(\alpha,\beta)*(\beta,\gamma)=(\alpha,\gamma)$. So one obtains
the natural bijections of sets
\begin{equation}\label{debut} 
  |K|_{flow}^0\iso  \liminj_{\square[n] \rightarrow K}
  (\{\widehat{0}<\widehat{1}\}^n)^0 \iso \liminj_{\square[n] \rightarrow
    K} \{\widehat{0},\widehat{1}\}^n \iso \liminj_{\square[n] \rightarrow
    K} \square[n]_0 \iso K_0.
\end{equation} 

\subsection{Small realization of a precubical set as a flow}
\label{secsmall}

The first two propositions will help the reader to understand the
differences between the realization functor $|-|_{flow}$ and the new
one $\gl(-)$ which is going to be constructed in this section.

\bp (e.g., \cite{MR99h:55031} Lemma~5.2.6) \label{cube} Let
$\mathcal{M}$ be a model category. Consider a pushout diagram of
$\mathcal{M}$
\[
\xymatrix{
A \fr{} \fd{} & B \fd{}\\
C \fr{} & \cocartesien D}
\]
such that the objects $A$, $B$ and $C$ are cofibrant and such that the
map $A \rightarrow C$ is a cofibration. Then $D$ is cofibrant and is
the homotopy colimit. In other terms, the commutative diagram above is
also a homotopy pushout diagram. \ep

\bp \label{ancien} The functor $|-|_{flow}:\square^{op}\set
\rightarrow \dtop$ is a left adjoint (and therefore is
colimit-preserving). Moreover, it satisfies the following properties:
\begin{itemize}
\item For every $n\geq 0$, there is a homotopy pushout diagram of
  flows
\[
\xymatrix{
\glob(\mathbf{S}^{n-1}) \fr{} \fd{}& |\de\square[n+1]|_{flow}\fd{}\\
\glob(\mathbf{D}^{n}) \fr{} & \hococartesien{|\square[n+1]|_{flow}}.}
\]
\item There exists an objectwise weak S-homotopy equivalence of
  cocubical flows \[|\square[*]|_{flow} \longrightarrow
  \{\widehat{0}<\widehat{1}\}^*\] (with always, by convention,
  $\{\widehat{0}<\widehat{1}\}^0=\{0\}$).
\end{itemize}
\ep 

\bpf One only has to prove the existence of the homotopy pushout
diagram. Equation~(\ref{debut}) implies $|\de\square[n+1]|_{flow}^0 =
\{\widehat{0},\widehat{1}\}^{n+1}$. By \cite{ccsprecub} Theorem~7.8,
there is a homotopy equivalence $\mathbf{S}^{n-1} \simeq
\P_{\widehat{0}\dots\widehat{0},\widehat{1}\dots\widehat{1}}|\de\square[n+1]|_{flow}$.
This yields a morphism of flows
\[t_n:\glob(\mathbf{S}^{n-1}) \rightarrow |\de\square[n+1]|_{flow}\]
defined by $t_n(\widehat{0})=\widehat{0}\dots\widehat{0}$,
$t_n(\widehat{1})=\widehat{1}\dots\widehat{1}$ and
\[\P_{\widehat{0},\widehat{1}}t_n :
\P_{\widehat{0},\widehat{1}}\glob(\mathbf{S}^{n-1}) = \mathbf{S}^{n-1}
\rightarrow
\P_{\widehat{0}\dots\widehat{0},\widehat{1}\dots\widehat{1}}\gl(\de\square[n+1])\] 
is a homotopy equivalence. Then consider the pushout diagram of flows: 
\[
\xymatrix{
\glob(\mathbf{S}^{n-1}) \fr{t_n} \fd{}& |\de\square[n+1]|_{flow}\fd{}\\
\glob(\mathbf{D}^{n}) \fr{} & \cocartesien{Z_{n+1}}.}
\]
By construction, one has the equality
$\P_{\alpha,\beta}|\square[n+1]|_{flow} =
\P_{\alpha,\beta}|\de\square[n+1]|_{flow}$ for every $(\alpha,\beta)
\neq (\widehat{0}\dots\widehat{0},\widehat{1}\dots\widehat{1})$ and
there is a pushout diagram of topological spaces
\[
\xymatrix{
  \mathbf{S}^{n-1} \ar@{->}[rr]^-{\P_{\widehat{0},\widehat{1}}t_n} \fd{} && \P_{\widehat{0}\dots\widehat{0},\widehat{1}\dots\widehat{1}}|\de\square[n+1]|_{flow} \fd{}\\
  \mathbf{D}^{n} \ar@{->}[rr] && \cocartesien
  \P_{\widehat{0}\dots\widehat{0},\widehat{1}\dots\widehat{1}}Z_{n+1}.
}
\]
Since the model category $\top$ is left proper, the map
$\mathbf{D}^{n} \rightarrow
\P_{\widehat{0}\dots\widehat{0},\widehat{1}\dots\widehat{1}}Z_{n+1}$
is a weak homotopy equivalence. So the flow $Z_{n+1}$ and $
|\square[n+1]|_{flow}$ are weakly S-homotopy equivalent. Since
$Z_{n+1}$ and $|\square[n+1]|_{flow}$ are both cofibrant and fibrant,
there is a S-homotopy equivalence $Z_{n+1}\simeq
|\square[n+1]|_{flow}$. The pushout diagram defining $Z_{n+1}$ is also
a homotopy pushout diagram by Proposition~\ref{cube}. Hence the
result.  \epf

\bp \label{ordrestate} Loopless flows satisfy the following two facts:
\begin{enumerate} 
\item If a flow $X$ is loopless, then the reflexive and transitive
  closure of the set
  \[\{(\alpha,\beta)\in X^0\p X^0\hbox{ such that
  }\P_{\alpha,\beta}X\neq\varnothing\}\] induces a partial ordering on
  $X^0$. 
\item The functor $X\mapsto X^0$ from flows to sets induces a functor
  from the full subcategory of loopless flows to that of partially
  ordered sets with strictly increasing maps.
\end{enumerate}
\ep

\bpf The first assertion is \cite{3eme} Lemma~4.2. The second
assertion is then clear. \epf 

\bth \label{small} There exists a colimit-preserving functor
$\gl:\square^{op}\set \rightarrow \dtop$ satisfying the following
properties: 
\begin{itemize}
\item For every $n\geq 0$, there is a pushout diagram of flows
\[
\xymatrix{
\glob(\mathbf{S}^{n-1}) \fr{} \fd{}& \gl(\de\square[n+1])\fd{}\\
\glob(\mathbf{D}^{n}) \fr{} & \cocartesien\gl(\square[n+1]).}
\]
\item There exists an objectwise weak S-homotopy equivalence of
  cocubical flows \[\gl(\square[*]) \longrightarrow
  \{\widehat{0}<\widehat{1}\}^*\] (with always, by convention,
  $\{\widehat{0}<\widehat{1}\}^0=\{0\}$). In particular with $n=0$,
  $\gl(\square[0])=\{0\}$.
\end{itemize}
\eth 

Note that by Proposition~\ref{cube}, the pushout diagram above is also
a homotopy pushout diagram.

\bpf Let us construct the restriction of the functor $\gl(-)$ to the
category of $n$-dimensional precubical sets $\square_n^{op}\set$ by
induction on $n\geq 0$. The functor $\gl(-)$ will satisfy the natural
isomorphism
\[\gl(K_{\leq n}) \iso \liminj_{\square[p] \rightarrow K_{\leq n}} \gl(\square[p])\] 
for every precubical set $K$. One will also prove by induction on
$n\geq 0$ that: 
\begin{itemize}
\item For any morphism $\delta$ of $\square_n$, the map
  $\gl(\square[\delta])$ is a relative $I^{gl}$-cell complex.
\item There exists a morphism of cocubical flows of dimension $n$ from
  $\gl([*])$ to $\{\widehat{0}<\widehat{1}\}^*$ which is an objectwise
  weak S-homotopy equivalence.
\item For all $0 \leq p\leq n$, the map $\gl(\de\square[p] \subset
  \square[p])$ is a cofibration. 
\end{itemize}

For $n=0$, let $\gl(\square[0])=\{0\}$. Note that this defines a
functor from $\square_0^{op}\set$ to $\dtop$ and that for any morphism
$\delta$ of $\square_0$, one has $\gl(\square[\delta])\in
\cell(I^{gl})$ since $\delta=\id_{[0]}$ is the only possibility.

Now suppose the construction done for $n\geq 0$. Consider the three
cocubical flows of dimension $n$ defined by $X([*]) =
\gl(\square[*])$, $Y([*]) = |\square[*]|_{flow}$ and $I([*]) =
\{\widehat{0}<\widehat{1}\}^*$ for all $0\leq * \leq n$. By induction
hypothesis, there exists a morphism of cocubical flows of dimension
$n$ from $X$ to $I$ which is an objectwise weak S-homotopy
equivalence. And Proposition~\ref{ancien} provides a morphism of
cocubical flows of dimension $n$ from $Y$ to $I$ which is an
objectwise weak S-homotopy equivalence. Any map from any cocubical
flow to $I$ is an objectwise fibration since for any $n\geq 0$, the
path space $\P I([n])$ is discrete. Finally, by induction hypothesis,
each map $\gl(\de\square[p] \subset \square[p])$ for $0 \leq p\leq n$
is a cofibration. And each map $|\de\square[p] \subset
\square[p]|_{flow}$ for $0 \leq p\leq n$ is a cofibration by
\cite{ccsprecub} Proposition~7.6.  Theorem~\ref{homotopy-natural-fini}
and Proposition~\ref{connexion} yield a natural S-homotopy equivalence
$\mu_{K_{\leq n}}:\gl(K_{\leq n}) \rightarrow |K_{\leq n}|_{flow}$.
The precubical set $\de\square[n+1]$ is of dimension $n$. So by
induction hypothesis, there exists a S-homotopy equivalence
\[\mu_{\de\square[n+1]}:\gl(\de\square[n+1])
\stackrel{\simeq}\longrightarrow |\de\square[n+1]|_{flow}.\] There is
also \[\gl(\de\square[n+1])^0\iso |\de\square[n+1]|_{flow}^0 \iso
\{\widehat{0},\widehat{1}\}^{n+1}\] by Equation~(\ref{debut}). The
continuous map
\[\P_{\widehat{0}\dots\widehat{0},\widehat{1}\dots\widehat{1}}\mu_{\de\square[n+1]}: 
\P_{\widehat{0}\dots\widehat{0},\widehat{1}\dots\widehat{1}}
\gl(\de\square[n+1]) \stackrel{\simeq}\longrightarrow
\P_{\widehat{0}\dots\widehat{0},\widehat{1}\dots\widehat{1}}
|\de\square[n+1]|_{flow}\] is a homotopy equivalence by \cite{model3}
Corollary~19.8. Using \cite{ccsprecub} Theorem~7.8, one deduces that
the topological spaces
$\P_{\widehat{0}\dots\widehat{0},\widehat{1}\dots\widehat{1}}
\gl(\de\square[n+1])$ and $\mathbf{S}^{n-1}$ are homotopy equivalent.
This yields a morphism of flows
\[s_n:\glob(\mathbf{S}^{n-1}) \rightarrow \gl(\de\square[n+1])\]
defined by
\begin{itemize}
\item $s_n(\widehat{0})=\widehat{0}\dots\widehat{0}$ with the
  identifications $\gl(\de\square[n+1])^0=\de\square[n+1]_0=\{\widehat{0},\widehat{1}\}^{n+1}$
\item $s_n(\widehat{1})=\widehat{1}\dots\widehat{1}$ with the
  identifications $\gl(\de\square[n+1])^0=\de\square[n+1]_0=\{\widehat{0},\widehat{1}\}^{n+1}$
\item $\P_{\widehat{0},\widehat{1}}s_n :
  \P_{\widehat{0},\widehat{1}}\glob(\mathbf{S}^{n-1}) =
  \mathbf{S}^{n-1} \rightarrow
  \P_{\widehat{0}\dots\widehat{0},\widehat{1}\dots\widehat{1}}\gl(\de\square[n+1])$
  is a homotopy equivalence.
\end{itemize}
The flow $\gl(\square[n+1])$ is then defined by the pushout diagram:
\begin{equation}\label{preleft}
\xymatrix{
\glob(\mathbf{S}^{n-1}) \fr{s_n} \fd{}& \gl(\de\square[n+1])\fd{}\\
\glob(\mathbf{D}^{n}) \fr{} & \cocartesien\gl(\square[n+1]).}
\end{equation}
Note that by construction, the map $\gl(\de\square[n+1] \subset
\square[n+1])$ is a cofibration. The $2(n+1)$ inclusions $\square[n]
\subset \de\square[n+1]$ yield the definition of the
$\gl(\delta_i^\alpha)$'s for all $\delta_i^\alpha:[n] \rightarrow
[n+1]$ with $1\leq i\leq n+1$ and $\alpha\in\{0,1\}$ as the composites
\[\gl(\delta_i^\alpha):\gl(\square[n]) \subset \gl(\de\square[n+1])
\longrightarrow \gl(\square[n+1]).\] Since the category
$\square_{n+1}$ is the quotient of the free category generated by the
$\delta_i^\alpha:[p-1] \rightarrow [p]$ for $1\leq p\leq n+1$ with
$1\leq i\leq p$ and $\alpha\in\{0,1\}$, by the cocubical relations,
one has to check the cocubical relation $\gl(\delta_j^\beta) \circ
\gl(\delta_i^\alpha) = \gl(\delta_i^\alpha) \circ
\gl(\delta_{j-1}^\beta)$ for $i<j$ for every map $\delta_j^\beta \circ
\delta_i^\alpha:[n-1]\rightarrow [n+1]$ and $\delta_i^\alpha \circ
\delta_{j-1}^\beta:[n-1]\rightarrow [n+1]$. By induction, each
morphism of flows $\gl(\delta_i^\alpha)$ is an inclusion of
$I^{gl}$-cell subcomplexes. The equality $\delta_j^\beta\circ
\delta_i^\alpha = \delta_i^\alpha \circ \delta_{j-1}^\beta$ implies
that the sources of $\gl(\delta_j^\beta) \circ \gl(\delta_i^\alpha)$
and $\gl(\delta_i^\alpha) \circ \gl(\delta_{j-1}^\beta)$ are the same
$I^{gl}$-cell subcomplex of $\gl(\square[n+1])$~\footnote{This
  argument is possible since every element of $\cell(I^{gl})$ is an
  (effective) monomorphism of flows by \cite{model3} Theorem~10.6.
  Indeed, a subcomplex of a relative $I^{gl}$-cell complex is then
  entirely determined by its set of cells by \cite{ref_model2}
  Proposition~10.6.10 and Proposition~10.6.11.}.  Hence the equality.
So the functor from $\square_{n}$ to $\dtop$ defined by $[p]\mapsto
\gl(\square[p])$ for $p\leq n$ is extended to a functor from
$\square_{n+1}$ to $\dtop$ defined by $[p]\mapsto \gl(\square[p])$ for
$p\leq n+1$. Let
\[\gl(K_{\leq n+1}) = \liminj_{\square[p] \rightarrow K_{\leq n+1}} \gl(\square[p])\]
for all precubical sets $K$. This construction extends the functor
$\gl:\square_{n}^{op}\set \rightarrow \dtop$ to a functor
$\gl:\square_{n+1}^{op}\set \rightarrow \dtop$.

It remains to prove that one has an objectwise weak S-homotopy
equivalence of cocubical flows of dimension $n+1$ from
$\gl(\square[*])$ to $\{\widehat{0}<\widehat{1}\}^*$ to complete the
induction and the proof. The map 
\[\gl(\delta_i^\alpha):\gl(\square[n]) \subset \gl(\de\square[n+1])
\longrightarrow \gl(\square[n+1])\]
induces a set map 
\[\gl(\delta_i^\alpha)^0:\gl(\square[n])^0 \subset \gl(\de\square[n+1])^0
\longrightarrow \gl(\square[n+1])^0.\] By Equation~(\ref{debut}) and
Proposition~\ref{ordrestate}, one obtains a strictly increasing set
map \[\gl(\delta_i^\alpha)^0: \{\widehat{0}<\widehat{1}\}^{n}
\rightarrow \{\widehat{0}<\widehat{1}\}^{n+1}.\] This yields a
morphism of cocubical flows of dimension $n+1$ from $\gl(\square[*])$
to $\{\widehat{0}<\widehat{1}\}^*$. It remains to prove that
$\gl(\square[n+1])$ is weakly S-homotopy equivalent to
$\{\widehat{0}<\widehat{1}\}^{n+1}$. By construction of
$\gl(\square[n+1])$, one has the equality
$\P_{\alpha,\beta}\gl(\square[n+1]) =
\P_{\alpha,\beta}\gl(\de\square[n+1])$ for every $(\alpha,\beta) \neq
(\widehat{0}\dots\widehat{0},\widehat{1}\dots\widehat{1})$ and there
is a pushout diagram of topological spaces
\[
\xymatrix{
  \mathbf{S}^{n-1} \ar@{->}[rr]^-{\P_{\widehat{0},\widehat{1}}s_n} \fd{} && \P_{\widehat{0}\dots\widehat{0},\widehat{1}\dots\widehat{1}}\gl(\de\square[n+1]) \fd{}\\
  \mathbf{D}^{n} \ar@{->}[rr] && \cocartesien
  \P_{\widehat{0}\dots\widehat{0},\widehat{1}\dots\widehat{1}}\gl(\square[n+1]).
}
\]
Since the map $\P_{\widehat{0},\widehat{1}}s_n$ is a weak homotopy
equivalence, and since the model category $\top$ is left proper, the
map $\mathbf{D}^{n}\rightarrow
\P_{\widehat{0}\dots\widehat{0},\widehat{1}\dots\widehat{1}}\gl(\square[n+1])$
is a weak homotopy equivalence.  \epf

\begin{cor} \label{small-cor} There exist a natural transformation
  $\mu:\gl(-) \rightarrow |-|_{flow}$ inducing for every precubical
  set $K$ a natural S-homotopy equivalence $\mu_K:\gl(K) \simeq
  |K|_{flow}$ and a natural transformation $\nu:|-|_{flow} \rightarrow
  \gl(-)$ inducing for every precubical set $K$ a natural S-homotopy
  equivalence $\nu_K:|K|_{flow} \simeq \gl(K)$ which is an inverse up
  to S-homotopy of $\mu_K$.
\end{cor}

\bpf Consider the three cocubical flows $X([*]) = \gl(\square[*])$,
$Y([*]) = |\square[*]|_{flow}$ and $I([*]) =
\{\widehat{0}<\widehat{1}\}^*$ for all $*\geq 0$.  Theorem~\ref{small}
and Proposition~\ref{ancien} yield objectwise weak S-homotopy
equivalences $X\rightarrow I$ and $Y\rightarrow I$. Since the path
space $\P I([n])$ is discrete for all $n\geq 0$, the two maps
$X\rightarrow I$ and $Y\rightarrow I$ are objectwise trivial
fibrations of cocubical flows. Then let us apply
Theorem~\ref{homotopy-natural} and let us notice that a simplicial
homotopy gives rise to a S-homotopy by Proposition~\ref{connexion}.
\epf

The following theorem characterizes realization functors:

\bth Let $X \rightarrow \{\widehat{0}<\widehat{1}\}^*$ be an
objectwise weak S-homotopy equivalence of cocubical flows. Assume that
for every $n \geq 0$, the map $\widehat{X}(\de\square[n]) \rightarrow
\widehat{X}(\square[n])$ is a cofibration. Then there exist natural
transformations $\mu: \gl \rightarrow \widehat{X}$ and $\nu:
\widehat{X}\rightarrow \gl$ inducing natural S-homotopy equivalences
which are inverse to each other up to S-homotopy. In particular, for
all $n\geq 0$, there is a homotopy pushout diagram of flows
\[
\xymatrix{
\glob(\mathbf{S}^{n-1}) \fr{} \fd{}& \widehat{X}(\de\square[n+1])\fd{}\\
\glob(\mathbf{D}^{n}) \fr{} & \hococartesien{\widehat{X}(\square[n+1])}}
\]
where the left-hand vertical map is the inclusion of flows
$\glob(\mathbf{S}^{n-1} \subset \mathbf{D}^{n})$ and the right-hand
vertical map $\widehat{X}(\de\square[n+1] \subset \square[n+1])$.
\eth

\bpf Since the path space $\P \{\widehat{0}<\widehat{1}\}^n$ of the
flow $\{\widehat{0}<\widehat{1}\}^n$ is discrete for all $n\geq 0$,
the map $X \rightarrow \{\widehat{0}<\widehat{1}\}^*$ is an objectwise
trivial fibration of flows. Then apply Theorem~\ref{homotopy-natural}
and Theorem~\ref{small} to obtain the natural transformations $\mu$
and $\nu$. One obtains the commutative diagram of flows
\[
\xymatrix{
\glob(\mathbf{S}^{n-1}) \fr{} \fd{}& \fd{}\gl(\de\square[n+1])\ar@{->}[rr]^-{\mu_{\de\square[n+1]}}&& \widehat{X}(\de\square[n+1])\fd{}\\
\glob(\mathbf{D}^{n}) \fr{} & \cocartesien
{\gl(\square[n+1])}\ar@{->}[rr]^-{\mu_{\square[n+1]}}&& {\widehat{X}(\square[n+1])},}
\]
and therefore the commutative diagram of flows 
\[
\xymatrix{
\glob(\mathbf{S}^{n-1}) \fr{} \fd{}&
\fd{}\gl(\de\square[n+1])\ar@{->}[rr]^-{\mu_{\de\square[n+1]}}&&
\widehat{X}(\de\square[n+1])\fd{}\ar@{=}[r] &
\widehat{X}(\de\square[n+1]) \fd{}\\
\glob(\mathbf{D}^{n}) \fr{} & \cocartesien
{\gl(\square[n+1])}\ar@/_20pt/[rrr]_-{\mu_{\square[n+1]}}\ar@{->}[rr]^-{\phi_n}&& \cocartesien
{T_n} \fr{\psi_n} & {\widehat{X}(\square[n+1])}.}
\]
The map $\phi_n$ is a weak S-homotopy equivalence since $\dtop$ is
left proper by \cite{2eme} Theorem~7.4. So by the two-out-of-three
property, the map $\psi_n$ is a weak S-homotopy equivalence as well.
Hence the homotopy pushout of flows.  \epf

\section{Realizing a precubical set as a small globular complex}
\label{glob}

\subsection{Globular complex} 

A globular complex is, like a $d$-space, a local pospace and a stream,
a topological space with an additional structure modeling time
irreversibility. We refer to \cite{model2} for further explanations
about the following list of definitions. The original definition of a
globular complex can be found in \cite{diCW} but this old definition
is slightly different and less tractable than the one of
\cite{model2}. So it will not be used.

\begin{figure}
\begin{center}
\includegraphics[width=7cm]{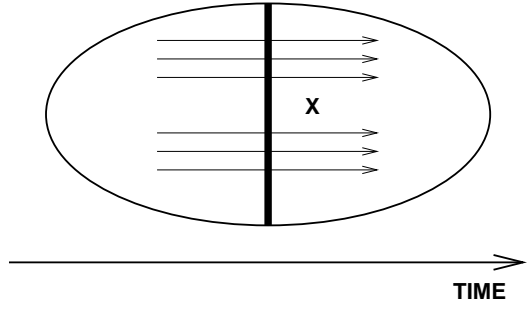}
\end{center}
\caption{Symbolic representation of
$\glob^{top}(X)$ for some topological space $X$} 
\label{exglob}
\end{figure}

A \textit{multipointed topological space} $(X,X^0)$ is a pair of
topological spaces such that $X^0$ is a discrete subspace of $X$.  A
morphism of multipointed topological spaces $f:(X,X^0)\longrightarrow
(Y,Y^0)$ is a continuous map $f:X\longrightarrow Y$ such that
$f(X^0)\subset Y^0$. The corresponding category is denoted by $\mtop$.
The category of multipointed spaces is cocomplete. Let $Z$ be a
topological space. The \textit{ (topological) globe of $Z$}, which is
denoted by $\glob^{top}(Z)$, is the multipointed space
$(|\glob^{top}(Z)|,\{\widehat{0},\widehat{1}\})$ where the topological
space $|\glob^{top}(Z)|$ is the quotient of
$\{\widehat{0},\widehat{1}\}\sqcup (Z\p[0,1])$ by the relations
$(z,0)=(z',0)=\widehat{0}$ and $(z,1)=(z',1)=\widehat{1}$ for any
$z,z'\in Z$ (cf. Figure~\ref{exglob}). In particular,
$\glob^{top}(\varnothing)$ is the multipointed space
$(\{\widehat{0},\widehat{1}\},\{\widehat{0},\widehat{1}\})$. If $Z$ is
a singleton, then the globe of $Z$ is denoted by $\vI^{top}$. Let
\[I^{gl,top}:=\{\glob^{top}(\mathbf{S}^{n-1})\longrightarrow
\glob^{top}(\mathbf{D}^{n}),n\geq 0\}.\] A \textit{globular
  precomplex} is a $\lambda$-sequence for some ordinal $\lambda$ of
multipointed topological spaces $X:\lambda\longrightarrow \mtop$ such
that $X \in\cell(I^{gl,top})$ and such that $X_0=(X^0,X^0)$ with $X^0$
a discrete space.  This $\lambda$-sequence is characterized by a
presentation ordinal $\lambda$, and for any $\beta<\lambda$ by an
integer $n_\beta\geq 0$ and an attaching map $\phi_\beta :
\glob^{top}(\mathbf{S}^{n_\beta-1}) \longrightarrow X_\beta$. The
family $(n_\beta,\phi_\beta)_{\beta<\lambda}$ is called the \textit{
  globular decomposition} of $X$. A morphism of globular precomplexes
$f:X\longrightarrow Y$ is a morphism of multipointed spaces still
denoted by $f$ from $\liminj X$ to $\liminj Y$. If $X$ is a globular
precomplex, then the \textit{underlying topological space} of the
multipointed space $\liminj X$ is denoted by $|X|$. Let $X$ be a
globular precomplex. A morphism of globular precomplexes
$\tau:\vI^{top}\longrightarrow X$ is \textit{non-decreasing} if there
exist $t_0=0<t_1<\dots<t_{n}=1$ such that:
\begin{enumerate}
\item  $\tau(t_i)\in X^0$ for all $0 \leq i \leq n$, 
\item $\tau(]t_i,t_{i+1}[)\subset \glob^{top}(\mathbf{D}^{n_{\beta_i}} 
\backslash \mathbf{S}^{n_{\beta_i}-1})$ for some $(n_{\beta_i},\phi_{\beta_i})$ 
of the globular decomposition of $X$, 
\item for $0\leq i<n$, there exists $z^i_\tau\in
  \mathbf{D}^{n_{\beta_i}}\backslash \mathbf{S}^{n_{\beta_i}-1}$ and a
  strictly increasing continuous map
  $\psi^i_\tau:[t_i,t_{i+1}]\longrightarrow [0,1]$ such that
  $\psi^i_\tau(t_i)=\widehat{0}$ and
  $\psi^i_\tau(t_{i+1})=\widehat{1}$ and for any $t\in
  [t_i,t_{i+1}]$, $\tau(t)=(z^i_\tau,\psi^i_\tau(t))$.
\end{enumerate}
In particular, the restriction $\tau\!\restriction_{]t_i,t_{i+1}[}$
of $\tau$ to $]t_i,t_{i+1}[$ is one-to-one. The set of
non-decreasing morphisms from $\vI^{top}$ to $X$ is denoted by
${\P}^{top}(X)$. A morphism of globular precomplexes
$f:X\longrightarrow Y$ is \textit{ non-decreasing} if the canonical
set map $\top([0,1],|X|)\longrightarrow \top([0,1],|Y|)$ induced by
composition by $f$ yields a set map ${\P}^{top}(X)\longrightarrow
{\P}^{top}(Y)$. In other terms, one has the commutative diagram of
sets
\[\xymatrix{
{\P}^{top}(X)\fr{}\fd{\subset}& {\P}^{top}(Y)\fd{\subset}\\
\top([0,1],|X|) \fr{} &\top([0,1],|Y|).}
\]
A \textit{globular complex} $X$ is a globular precomplex such that the
attaching maps $\phi_\beta$ are non-decreasing. A morphism of globular
complexes is a morphism of globular precomplexes which is
non-decreasing. The category of globular complexes together with the
morphisms of globular complexes as defined above is denoted by
$\gltop$.

\subsection{S-homotopy equivalence of globular complex} 

Let $X$ and $U$ be two globular complexes. Let $\tgltop(X,U)$ be the
set $\gltop(X,U)$ equipped with the Kelleyfication of the compact-open
topology. Let $f,g:X\rightrightarrows U$ be two morphisms of globular
complexes.  Then a \textit{S-homotopy} is a continuous map $H:[0,1]
\rightarrow \tgltop(X,U)$ with $H_0=f$ and $H_1=g$. This situation is
denoted by $f\sim_S g$.  The S-homotopy relation defines a congruence
on the category $\gltop$.  If there exists a map $f':U\rightarrow X$
with $f\circ f'\sim_S \id_U$ and $f'\circ f\sim_S \id_X$, then $f$ is
called a \textit{S-homotopy equivalence}. The class of S-homotopy
equivalences of globular complexes is denoted by $\mathcal{SH}$.

Since the S-homotopy relation of globular complexes is associated with
a cylinder functor (\cite{model2} Corollary~II.4.9), there is an
isomorphism of categories \[\gltop[\mathcal{SH}^{-1}] \iso
\gltop/\!\sim_S\] between the localization of the category of globular
complexes by the S-homotopy equivalences and the quotient of the
category of globular complexes by S-homotopy (see the proof of
\cite{model2} Theorem~V.4.1 and also \cite{model3} Theorem~4.7).

\subsection{Realizing a globular complex as a flow}

By \cite{model2} Theorem~III.3.1, there exists a unique functor
$cat:\gltop\longrightarrow\hda \subset \dtop$ such that
\begin{enumerate}
\item if $X=X^0$ is a discrete globular complex, then $cat(X)$ is the
  flow $X^0$
\item if $Z=\mathbf{S}^{n-1}$ or $Z=\mathbf{D}^{n}$ for some integer $n\geq 0$, 
then $cat(\glob^{top}(Z))=\glob(Z)$, 
\item for any globular complex $X$ with globular decomposition
  $(n_\beta,\phi_\beta)_{\beta<\lambda}$, for any limit ordinal
  $\beta\leq\lambda$, the canonical morphism of flows
\[\liminj_{\alpha<\beta} cat(X_\alpha)\longrightarrow
cat(X_\beta)\] is an isomorphism of flows, 
\item for any globular complex $X$ with globular decomposition
  $(n_\beta,\phi_\beta)_{\beta<\lambda}$, for any $\beta<\lambda$, one
  has the pushout of flows
\[\xymatrix{\glob(\mathbf{S}^{n_\beta-1})\fr{cat(\phi_\beta)}\fd{}& cat(X_\beta)\fd{}\\
  \glob(\mathbf{D}^{n_\beta})\fr{} & cat(X_{\beta+1}).\cocartesien}\]
\end{enumerate}
The properties of the functor $cat$ used in this paper are summarized
in the statement below:

\bth \label{rappelcat} One has:
\begin{itemize}
\item The functor $cat$ induces for each pair $(X,U)$ of globular
  complexes a surjective set map $\gltop(X,U) \rightarrow
  \dtop(cat(X),cat(U))$ (\cite{model2} Corollary~IV.3.15).
\item For each flow $X\in\hda$, there exists a globular complex
  $X^{top}$ with $cat(X^{top})=X$ (\cite{model2} Theorem~V.4.1 and
  \cite{4eme} Theorem~6.1).
\item The functor $cat:\gltop \rightarrow \dtop$ induces a category
  equivalence \[\gltop[\mathcal{SH}^{-1}] \simeq \ho(\dtop)\] between
  the localization of $\gltop$ by the S-homotopy equivalences and the
  homotopy category of flows (\cite{model2} Theorem~V.4.2).
\item There exists a unique functor $|-|:\gltop[\mathcal{SH}^{-1}]
  \rightarrow \ho(\top)$ such that the following diagram of
  categories is commutative:
  \[
\xymatrix{
\gltop \fr{|-|} \fd{\gamma_\gltop} & \top \fd{\gamma_\top}\\
\gltop[\mathcal{SH}^{-1}] \fr{|-|} & \ho(\top)}
\]
where $\gamma_\gltop:\gltop \rightarrow \gltop[\mathcal{SH}^{-1}]$ is
the canonical functor from the category of globular complexes to its
localization by the S-homotopy equivalences (\cite{model2}
Corollary~VII.2.3).
\end{itemize}
\eth 

\subsection{Realizing a precubical set as a small globular complex}

\bp (Unknown reference) \label{compose} Let $\C$ be a cocomplete
category. Consider a commutative diagram of $\C$:
\[
\xymatrix{
A \fr{} \fd{} & B \fr{} \fd{} & C \fd{} \\
D \fr{}  & \cocartesien E \fr{} & F}
\]
such the square $ABDE$ is a pushout diagram. Then the square $ACDF$ is
a pushout diagram if and only if the square $BCEF$ is a pushout
diagram.  \ep

\bpf Well-known: 1) if $BCEF$ is a pushout diagram, from a cone
$(D\leftarrow A\rightarrow C)\stackrel{\bullet}\rightarrow Z$, one
deduces a cone $(E\leftarrow B\rightarrow
C)\stackrel{\bullet}\rightarrow Z$ and a map $F\rightarrow Z$; 2) if
$ACDF$ is a pushout diagram, from a cone $(E\leftarrow B\rightarrow
C)\stackrel{\bullet}\rightarrow Z$, one deduces a cone $(D\leftarrow
A\rightarrow C)\stackrel{\bullet}\rightarrow Z$, and a map
$F\rightarrow Z$; the composite $E\rightarrow F \rightarrow Z$ is
equal to the map $E\rightarrow Z$ of the cone $(E\leftarrow
B\rightarrow C)\stackrel{\bullet}\rightarrow Z$ since the two maps
satisfy the universal property relative to the pushout square $ABDE$;
so $F\rightarrow Z$ is a universal solution for the cone $(E\leftarrow
B\rightarrow C)\stackrel{\bullet}\rightarrow Z$.  \epf

\bth \label{reaglob} There exists a functor
$\gl^{top}:\square^{op}\set \rightarrow \gltop$ such that for every
precubical set $K$, there is a natural isomorphism of flows
$cat(\gl^{top}(K))\iso \gl(K)$. \eth

The functor $\gl^{top}(-)$ which is going to be constructed
\textit{essentially} coincides with the functor from precubical sets
to globular complexes constructed in \cite{diCW}. Essentially means
not exactly. Indeed, the old definition of globular complex given in
\cite{diCW} and the new one given in \cite{model2} are not exactly the
same. For example, with the new definition, an execution path is
locally strictly increasing: see the remark in \cite{model2} between
Definition~II.2.14 and Definition~II.2.15. Another difference:
Hausdorff spaces are used in \cite{diCW}. Weak Hausdorff spaces are
used here and in \cite{model2}.  Moreover, the construction given in
the following proof is more tractable than the construction given in
\cite{diCW} thanks to the use of the cocomplete category of
multipointed topological spaces. Note that the functor $cat$ is not
colimit-preserving. It only preserves globular decompositions of
globular complexes. So the proof is a little bit more complicated than
expected. Intuitively, the construction of $\gl^{top}(K)$ consists of
replacing each globe $\glob(\mathbf{D}^n)$ of $\gl(K)$ by a
topological globe $\glob^{top}(\mathbf{D}^n)$.

\bpf[Proof of Theorem~\ref{reaglob}] First of all, let us construct
the restriction of the functor $\gl^{top}(-)$ to $\square_n^{op}\set$
and let us prove the existence of a natural isomorphism
$cat(\gl^{top}(K_{\leq n}))\iso \gl (K_{\leq n})$ by induction on
$n\geq 0$. The functor $\gl^{top}(-)$ will satisfy the natural
isomorphism
\[\gl^{top}(K_{\leq n}) \iso \liminj_{\square[p] \rightarrow K_{\leq n}} \gl^{top}(\square[p])\] 
for every precubical set $K$, where the colimit is taken in the
category of multipointed topological spaces $\mtop$. So viewed as a
functor from $\square_n^{op}\set$ to $\mtop$, the functor $\gl^{top}(-)$
is a left adjoint. We will also prove by induction on $n$ that for any
morphism $\delta$ of $\square_n$, the morphism of globular complexes
$\gl^{top}(\square(\delta))$ is an element of $\cell(I^{gl,top})$.

For $n=0$, let $\gl^{top}(K_{\leq 0}) = K_0$. We have done since
$cat(\gl^{top}(K_{\leq 0})) = K_0$. Note this defines a functor from
$\square_0^{op}\set$ to $\gltop$ which is colimit-preserving. Note
also that for any morphism $\delta$ of $\square_0$, one has
$\gl^{top}(\square[\delta])\in \cell(I^{gl,top})$, $\delta=\id_{[0]}$
being the only possibility.

Now suppose the construction done for $n\geq 0$.  The precubical set
$\de\square[n+1]$ is of dimension $n$. So the globular complex
$\gl^{top}(\de\square[n+1])$ is already defined by induction
hypothesis and one has the isomorphism of flows
$cat(\gl^{top}(\de\square[n+1]))\iso \gl(\de\square[n+1])$. Since the
set map
\[\gltop(\glob^{top}(\mathbf{S}^{n-1}),\gl^{top}(\de\square[n+1]))
\rightarrow
\dtop(\glob(\mathbf{S}^{n-1}),\gl(\de\square[n+1]))\] is onto
by Theorem~\ref{rappelcat}, there exists a morphism of globular
complexes
\[s_n^{top}:\glob^{top}(\mathbf{S}^{n-1}) \rightarrow
\gl^{top}(\de\square[n+1])\] with $cat(s_n^{top})=s_n$, $s_n$ being
the map defined in the proof of Theorem~\ref{small}.  Let
$\gl^{top}(\square[n+1])$ be the multipointed topological space
defined by the pushout diagram of multipointed topological spaces
\begin{equation}\label{left}
\xymatrix{
\glob^{top}(\mathbf{S}^{n-1})\fr{s_n^{top}} \fd{} & \gl^{top}(\de\square[n+1])
\fd{} \\
\glob^{top}(\mathbf{D}^{n})\fr{} & \cocartesien\gl^{top}(\square[n+1]).}
\end{equation}
The globular decomposition of the multipointed space
$\gl^{top}(\square[n+1])$ is obtained by considering the globular
decomposition of the globular complex $\gl^{top}(\de\square[n+1])$ and
by adding the globular cell $\glob^{top}(\mathbf{S}^{n-1})\subset
\glob^{top}(\mathbf{D}^{n})$ with the attaching map $s_n^{top}$. So
the multipointed space $\gl^{top}(\square[n+1])$ is a globular
complex.

The $2(n+1)$ inclusions $\square[n] \subset \de\square[n+1]$ yield the
definition of the $\gl^{top}(\delta_i^\alpha)$'s for all
$\delta_i^\alpha:[n] \rightarrow [n+1]$ with $1\leq i\leq n+1$ and
$\alpha\in\{0,1\}$ as the composites:
\[\gl^{top}(\delta_i^\alpha):\gl^{top}(\square[n]) \subset \gl^{top}(\de\square[n+1])
\longrightarrow \gl^{top}(\square[n+1]).\] Since the category
$\square_{n+1}$ is the quotient of the free category generated by the
$\delta_i^\alpha:[p-1] \rightarrow [p]$ for $1\leq p\leq n+1$ with
$1\leq i\leq p$ and $\alpha\in\{0,1\}$, by the cocubical relations,
one has to check the cocubical relation $\gl^{top}(\delta_j^\beta)
\circ \gl^{top}(\delta_i^\alpha) = \gl^{top}(\delta_i^\alpha) \circ
\gl^{top}(\delta_{j-1}^\beta)$ for $i<j$ for every map $\delta_j^\beta
\circ \delta_i^\alpha:[n-1]\rightarrow [n+1]$ and $\delta_i^\alpha
\circ \delta_{j-1}^\beta:[n-1]\rightarrow [n+1]$. By induction, each
morphism of globular complexes $\gl^{top}(\delta_i^\alpha)$ is an
inclusion of $I^{gl,top}$-cell subcomplexes. The equality
$\delta_j^\beta\circ \delta_i^\alpha = \delta_i^\alpha \circ
\delta_{j-1}^\beta$ implies that the sources of
$\gl^{top}(\delta_j^\beta) \circ \gl^{top}(\delta_i^\alpha)$ and
$\gl^{top}(\delta_i^\alpha) \circ \gl^{top}(\delta_{j-1}^\beta)$ are
the same $I^{gl,top}$-cell subcomplex of
$\gl^{top}(\square[n+1])$~\footnote{This argument is possible since
  every element of $\cell(I^{gl,top})$ is an (effective) monomorphism
  of multipointed topological spaces by \cite{4eme} Theorem~8.2.
  Indeed, since $\mtop$ is cocomplete, a subcomplex of a relative
  $I^{gl,top}$-cell complex is then entirely determined by its set of
  cells by \cite{ref_model2} Proposition~10.6.10 and
  Proposition~10.6.11.}.  Hence the equality.

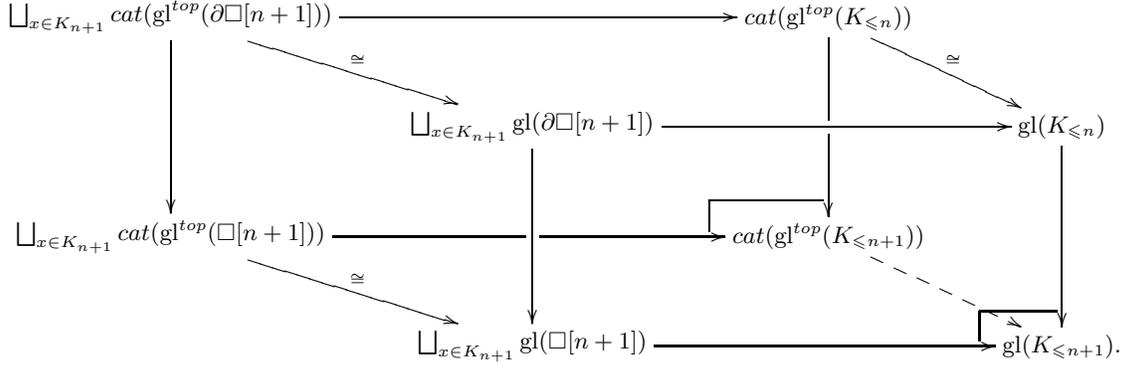
\begin{figure}
{\footnotesize
\[
\xymatrix{ \bigsqcup_{x\in K_{n+1}} cat(\gl^{top}(\de\square[n+1]))
  \ar@{->}[rd]^-{\iso} \ar@{->}[dd] \ar@{->}[rr] && cat(\gl^{top}(K_{\leq
    n}))
  \ar@{->}[rd]^-{\iso} \ar@{->}'[d][dd]  &\\
  & \bigsqcup_{x\in K_{n+1}}
  \gl(\de\square[n+1]) \ar@{->}[rr] \ar@{->}[dd] && \gl(K_{\leq n}) \ar@{->}[dd] \\
  \bigsqcup_{x\in K_{n+1}}
  cat(\gl^{top}(\square[n+1])) \ar@{->}[rd]^-{\iso} \ar@{->}'[r][rr] &&
  \cocartesien cat(\gl^{top}(K_{\leq n+1})) \ar@{-->}[rd]^-{}& \\
  & \bigsqcup_{x\in K_{n+1}} \gl(\square[n+1]) \ar@{->}[rr] &&
  \cocartesien \gl(K_{\leq n+1}).}
\]
} 
\caption{Isomorphism $cat(\gl^{top}(K_{\leq n+1})) \iso \gl(K_{\leq
    n+1})$.}
\label{cube3}
\end{figure}

So the functor from $\square_{n}$ to $\gltop$ defined by $[p]\mapsto
\gl^{top}(\square[p])$ for $p\leq n$ is extended to a functor from
$\square_{n+1}$ to $\mtop$ defined by $[p]\mapsto
\gl^{top}(\square[p])$ for $p\leq n+1$. Let
\[\gl^{top}(K_{\leq n+1}) = \liminj_{\square[p] \rightarrow K_{\leq n+1}} \gl^{top}(\square[p])\]
for all precubical sets $K$. This construction extends the functor
$\gl^{top}:\square_{n}^{op}\set \rightarrow \gltop$ to a functor
$\gl^{top}:\square_{n+1}^{op}\set \rightarrow \mtop$ which is still
colimit-preserving since it is still a left adjoint. So one obtains
the commutative diagram of multipointed spaces
\[
\xymatrix{
\bigsqcup_{x\in K_{n+1}} \glob^{top}(\mathbf{S}^{n-1})\fr{\bigsqcup s_n^{top}} \fd{} &
\bigsqcup_{x\in K_{n+1}} \gl^{top}(\de\square[n+1])  \fr{} \fd{} &
\gl^{top}(K_{\leq n})
\fd{} \\
\bigsqcup_{x\in K_{n+1}} \glob^{top}(\mathbf{D}^{n})\fr{} &
\bigsqcup_{x\in K_{n+1}} \gl^{top}(\square[n+1]) \cocartesien\fr{} &
\cocartesien\gl^{top}(K_{\leq n+1}).}
\] 
The left-hand square is a pushout by definition of
$\gl^{top}(\square[n+1])$. The right-hand square is a pushout since
$\gl^{top}:\square^{op}\set \rightarrow \mtop$ is colimit-preserving
and since for every precubical set $K$, there is a pushout diagram of
sets
\[
\xymatrix{ \bigsqcup_{x\in K_{n+1}} \de\square[n+1] \fr{} \fd{} &
  K_{\leq n}
  \fd{} \\
  \bigsqcup_{x\in K_{n+1}} \square[n+1] \fr{} & \cocartesien K_{\leq
    n+1}}
\] 
where the sum is over $x\in K_{n+1}=\square^{op}\set(\square[n+1],K)$
and where the corresponding map $\de\square[n+1] \rightarrow K_{\leq
  n}$ is the composite $\de\square[n+1] \subset \square[n+1]
\stackrel{x}\rightarrow K_{\leq n}$. Since pushout diagrams compose by
Proposition~\ref{compose}, one obtains the pushout diagram of
multipointed spaces
\begin{equation}\label{fabrication_glob}
\xymatrix{
\bigsqcup_{x\in K_{n+1}} \glob^{top}(\mathbf{S}^{n-1})\fr{} \fd{} &
\gl^{top}(K_{\leq n})
\fd{} \\
\bigsqcup_{x\in K_{n+1}} \glob^{top}(\mathbf{D}^{n})\fr{}  &
\cocartesien\gl^{top}(K_{\leq n+1}),}
\end{equation}
an then, by construction of the functor $cat:\gltop\rightarrow \dtop$,
the pushout diagram of flows
\[
\xymatrix{
\bigsqcup_{x\in K_{n+1}} \glob(\mathbf{S}^{n-1})\fr{} \fd{} &
cat(\gl^{top}(K_{\leq n}))
\fd{} \\
\bigsqcup_{x\in K_{n+1}} \glob(\mathbf{D}^{n})\fr{}  &
\cocartesien cat(\gl^{top}(K_{\leq n+1})).}
\]
Diagram~(\ref{fabrication_glob}) yields a globular decomposition for
the multipointed space $\gl^{top}(K_{\leq n+1})$, proving that the
functor $\gl^{top}(-)$ is actually a functor from
$\square_{n+1}^{op}\set$ to $\gltop$. By construction of the functor
$cat:\gltop\rightarrow \dtop$, there is a pushout diagram of flows
\[
\xymatrix{
\bigsqcup\limits_{x\in K_{n+1}} \glob(\mathbf{S}^{n-1})\fr{\bigsqcup\limits s_n} \fd{} &
\bigsqcup\limits_{x\in K_{n+1}} cat(\gl^{top}(\de\square[n+1]))   \fd{} \\
\bigsqcup\limits_{x\in K_{n+1}} \glob(\mathbf{D}^{n})\fr{} & \cocartesien
\bigsqcup\limits_{x\in K_{n+1}} cat(\gl^{top}(\square[n+1])).}
\]
So by Proposition~\ref{compose}, one obtains the pushout diagram of
flows~\footnote{Let us repeat that the functor $cat$ is not
  colimit-preserving. So the use of Proposition~\ref{compose} seems to
  be necessary to obtain Diagram~(\ref{back}).}:
\begin{equation} \label{back}
\xymatrix{
\bigsqcup_{x\in K_{n+1}} cat(\gl^{top}(\de\square[n+1]))  \fr{} \fd{} &
cat(\gl^{top}(K_{\leq n}))
\fd{} \\
\bigsqcup_{x\in K_{n+1}} cat(\gl^{top}(\square[n+1])) \fr{} &
\cocartesien cat(\gl^{top}(K_{\leq n+1})).}
\end{equation} 
The diagram of solid arrows of Figure~\ref{cube3} is commutative for
the following reasons:
\begin{itemize}
\item The back face is commutative and is a pushout diagram of flows
  by Diagram~(\ref{back}).
\item The front face is commutative and is a pushout diagram of flows
  since the functor $\gl:\square^{op}\set\rightarrow \dtop$ is
  colimit-preserving.
\item Apply the functor $cat$ to Diagram~(\ref{left}). One obtains the
  pushout diagram of flows
\[
\xymatrix{
\glob(\mathbf{S}^{n-1})\ar@{->}[rr]^-{cat(s_n^{top})} \fd{} &&
cat(\gl^{top}(\de\square[n+1]))
\fd{} \\
\glob(\mathbf{D}^{n})\ar@{->}[rr]  &&
\cocartesien cat(\gl^{top}(\square[n+1])).}
\]
Diagram ~(\ref{preleft}) and the equality $cat(s_n^{top})=s_n$ imply
the commutativity of the left-hand face.
\item Finally, the top face is commutative since there is a natural
  isomorphism \[cat(\gl^{top}(K))\iso \gl(K)\] for all precubical sets
  $K$ of dimension $n$ by hypothesis.
\end{itemize}
Hence the existence of an isomorphism of flows
\[cat(\gl^{top}(K_{\leq n+1})) \iso \gl(K_{\leq n+1})\] for every
precubical set $K$. The isomorphism is natural for the following
reasons:
\begin{itemize}
\item The map $\bigsqcup_{x\in K_{n+1}}
  cat(\gl^{top}(\de\square[n+1])) \rightarrow cat(\gl^{top}(K_{\leq
    n}))$ is natural with respect to $K$ since it is the image by the
  functor $cat\circ \gl^{top}(-)$ of the natural map of precubical
  sets $i(K,n):\bigsqcup_{x\in K_{n+1}} \de\square[n+1] \rightarrow K_{\leq
    n}$.
\item The map $\bigsqcup_{x\in K_{n+1}} \gl(\de\square[n+1])
  \rightarrow \gl(K_{\leq n})$ is natural with respect to $K$ since it
  is the image by the functor $\gl(-)$ of the natural map of precubical
  sets $i(K,n)$.
\item There is a natural isomorphism $cat(\gl^{top}(L))\iso \gl(L)$
  with respect to $L$ for every $n$-dimensional precubical set $L$ by
  induction hypothesis. Apply this fact for $L=K_{\leq n}$ and
  $L=\bigsqcup_{x\in K_{n+1}} \de\square[n+1]$. 
\item Morphisms of precubical sets of the form $\bigsqcup_{x\in
    K_{n+1}} A \rightarrow \bigsqcup_{x\in K_{n+1}} B$ are natural
  with respect to $K$ for every morphism of precubical sets
  $A\rightarrow B$.
\end{itemize}

The induction is now complete. One has $cat(\gl^{top}(K)) \iso \liminj
cat(\gl^{top}(K_{\leq n}))$ by definition of the functor $cat$. And
one has $\gl(K) \iso \liminj \gl(K_{\leq n})$ since the functor
$\gl(-)$ is colimit-preserving. Hence a natural isomorphism of flows
$cat(\gl^{top}(K)) \iso\gl(K)$.  \epf

Note that the functor $\gl^{top}:\square^{op}\set\rightarrow \gltop$
is not unique. It is entirely characterized up to isomorphism of
functors by the non-canonical choice of the maps
$s_n:\glob(\mathbf{S}^{n-1}) \rightarrow \gl(\de\square[n+1])$ and of
the maps $s_n^{top}:\glob^{top}(\mathbf{S}^{n-1}) \rightarrow
\gl^{top}(\de\square[n+1])$ for all $n\geq 0$. Let $\gamma_\dtop:\dtop
\rightarrow \ho(\dtop)$ be the canonical functor from the category of
flows to its homotopy category. Let us denote by
\[\overline{cat}:\gltop[\mathcal{SH}^{-1}] \simeq \ho(\top):\overline{cat}^{-1}\] the
equivalence of categories between the globular complexes up to
S-homotopy and the homotopy category of flows (see
Theorem~\ref{rappelcat}).

\bth \label{unique} The functor $\gamma_\dtop\circ
|-|_{flow}:\square^{op}\set \rightarrow \ho(\dtop)$ factors up to an
isomorphism of functors as a composite \[\square^{op}\set
\stackrel{\hogl^{top}}\longrightarrow \gltop/\!\sim_S \longrightarrow
\ho(\dtop).\] The functor $\hogl^{top}: \square^{op}\set \rightarrow
\gltop/\sim_S$ is unique up to isomorphism of functors.  \eth

In other terms, the functor $\gl^{top}(-)$ constructed in
Theorem~\ref{reaglob} is unique up to a natural S-homotopy of globular
complexes.

\bpf Since there is an isomorphism of categories
$\gltop[\mathcal{SH}^{-1}] \iso \gltop/\!\sim_S$, let us identify the
two categories.  Let \[\boxed{\hogl^{top}=\gamma_\gltop\circ
  \gl^{top}}.\] Then one obtains the isomorphisms of functors
$\overline{cat}\circ \hogl^{top}=\gamma_\dtop\circ cat \circ \gl^{top}
\iso \gamma_\dtop\circ |-|_{flow}$ by Theorem~\ref{reaglob} and
Theorem~\ref{small}. Hence the existence. Take two functors
$\hogl_1^{top}: \square^{op}\set \rightarrow \gltop/\!\sim_S$ and
$\hogl_2^{top}: \square^{op}\set \rightarrow \gltop/\!\sim_S$
satisfying the condition of the theorem.  Then $\overline{cat}\circ
\hogl_1^{top} = \overline{cat}\circ \hogl_2^{top} = \gamma_\dtop\circ
|-|_{flow}$. So one has the isomorphisms of functors $\hogl_1^{top}
\iso \overline{cat}^{-1}\circ (\overline{cat}\circ \hogl_1^{top}) \iso
\overline{cat}^{-1}\circ (\overline{cat}\circ \hogl_2^{top}) \iso
\hogl_2^{top}$. \epf

\section{Globular and cubical underlying homotopy type}
\label{under}

\subsection{Definition of the globular and cubical underlying homotopy type}

Let $\square \rightarrow \top$ be the functor defined on objects by
the mapping $[n]\mapsto [0,1]^n$ and on morphisms by the mapping $
\delta_i^\alpha \mapsto \lp (\epsilon_1, \dots, \epsilon_{n-1})
\mapsto (\epsilon_1, \dots, \epsilon_{i-1}, \alpha, \epsilon_i, \dots,
\epsilon_{n-1})\rp$. The functor $| - |_{space}: \square^{op}\set
\rightarrow \top$ is then defined by \[\boxed{|K|_{space}:=
  \liminj_{\square[n] \rightarrow K} [0,1]^n}.\] It is a left adjoint.
So it commutes with all small colimits. The purpose of this section is
the comparison of this functor with the \textit{underling homotopy
  type functor} defined by the composite \cite{model2}:
\[\boxed{
\xymatrix@1{
\Omega:\dtop \fr{\gamma_\dtop} & \ho(\dtop) \fr{\overline{cat}^{-1}} & 
\gltop[\mathcal{SH}^{-1}] \fr{|-|} & \ho(\top).}}\]

\subsection{Comparison of the two functors}

\bth \label{app} For every precubical set $K$, there is a natural
isomorphism of homotopy types $\gamma_\top(|K|_{space}) \iso
\Omega(|K|_{flow})$. \eth

\bpf Consider the three cocubical topological spaces
\begin{itemize}
\item $X([*])=[0,1]^*$ for all $*\geq 0$
\item $Y([*])=|\gl^{top}(\square[*])|$ for all $*\geq 0$
\item $I([*])=\{0\}$ for all $*\geq 0$.
\end{itemize} 
There exist a unique map $X\rightarrow I$ and a unique map
$Y\rightarrow I$ which are both objectwise weak homotopy equivalences
and objectwise fibrations.  The cocubical object $X$ satisfies the
hypotheses of Theorem~ \ref{homotopy-natural} in an obvious way. The
proof of Theorem~\ref{reaglob} implies the pushout diagram of spaces
\[
\xymatrix{
|\glob^{top}(\mathbf{S}^{n-1})|\fr{s_n^{top}} \fd{} & |\gl^{top}(\de\square[n+1])|
\fd{} \\
|\glob^{top}(\mathbf{D}^{n})|\fr{} & \cocartesien{|\gl^{top}(\square[n+1])|}.}
\] 
for every $n\geq 0$. Since the continuous map
$|\glob^{top}(\mathbf{S}^{n-1})| \rightarrow
|\glob^{top}(\mathbf{D}^{n})|$ is a cofibration of topological spaces
(see the proof of \cite{4eme} Theorem~8.2), the map
$|\gl^{top}(\de\square[n+1])| \rightarrow |\gl^{top}(\square[n+1])|$
is a cofibration as well. So the cocubical space $Y$ satisfies the
hypotheses of Theorem~\ref{homotopy-natural} as well. Hence there
exists a natural homotopy equivalence $|K|_{space} \simeq
|\gl^{top}(K)|$ with respect to $K$.  So there exists a natural
isomorphism of homotopy types $\gamma_\top(|K|_{space}) \iso
\gamma_\top(|\gl^{top}(K)|)$.  The proof is complete after the
following sequence of natural isomorphisms:
\begin{align*}
  &\gamma_\top(|K|_{space}) &\\ & \iso \gamma_\top(|\gl^{top}(K)|) &
  \hbox{ by the result above}\\ &\iso |\gamma_\gltop(\gl^{top}(K))| &
  \hbox{ by Theorem~\ref{rappelcat}}\\ &\iso |\overline{cat}^{-1}\circ
  \overline{cat}\circ \gamma_\gltop(\gl^{top}(K))| & \hbox{ since
    $\overline{cat}^{-1}\circ
    \overline{cat}\iso\id_{\gltop[\mathcal{SH}^{-1}]}$}\\ & \iso
  |\overline{cat}^{-1}\circ\gamma_\dtop(cat(\gl^{top}(K)))|& \hbox{ by
    Theorem~\ref{rappelcat}}\\ & \iso
  |\overline{cat}^{-1}\circ\gamma_\dtop(|K|_{flow})| & \hbox{ by
    Theorem~\ref{reaglob}}\\ & = \Omega(|K|_{flow})& \hbox{ by
    definition of $\Omega$.}\\
\end{align*}
\epf

\end{document}